\documentclass[11pt,amsfonts, fullpage]{amsart}

\usepackage{fullpage}

\def\endproof{\qed \medskip}
\def\blacksquare{\hbox to .60em {\vrule width .60em height .60em}}

\newtheorem{theorem}{Theorem}[section]

\newtheorem{corollary}[theorem]{Corollary}

\newtheorem{lemma}[theorem]{Lemma}

\newtheorem{proposition}[theorem]{Proposition}

\newtheorem{remark}[theorem]{Remark}

\begin{document}

\title[]{On Boundary Value Problems for Einstein Metrics}

\author[]{Michael T. Anderson}

\thanks{Partially supported by NSF Grant DMS 0604735}

\maketitle

\abstract
On any given compact manifold $M^{n+1}$ with boundary $\partial M$, it is 
proved that the moduli space ${\mathcal E}$ of Einstein metrics on $M$, if 
non-empty, is a smooth, infinite dimensional Banach manifold, at least when 
$\pi_{1}(M, \partial M) = 0$. Thus, the Einstein moduli space is unobstructed. 
The usual Dirichlet and Neumann boundary maps to data on $\partial M$ are 
smooth, but not Fredholm. Instead, one has natural mixed boundary-value 
problems which give Fredholm boundary maps. 

   These results also hold for manifolds with compact boundary which have 
a finite number of locally asymtotically flat ends, as well as for the 
Einstein equations coupled to many other fields. 
\endabstract

\setcounter{section}{0}

\section{Introduction.}
\setcounter{equation}{0}

 Let $M = M^{n+1}$ be a compact $(n+1)$-dimensional manifold with 
boundary $\partial M$, $n \geq 2$. In this paper, we consider the structure 
of the space of Einstein metrics on $(M, \partial M)$, i.e.~metrics $g$ on 
$\bar M = M \cup \partial M$ satisfying the Einstein equations 
\begin{equation} \label{e.1.1}
Ric_{g} = \lambda g.
\end{equation}
Here $\lambda$ is a fixed constant, equal to $\frac{s}{n+1}$, where $s$ is 
the scalar curvature. It is natural to consider boundary value problems 
for the equations (1.1). For example, the Dirichlet problem asks: given a 
(smooth) Riemannian metric $\gamma$ on $\partial M$, determine whether there 
exists a Riemannian metric $g$ on $\bar M$, which satisfies the Einstein equations 
(1.1) with the boundary condition
\begin{equation} \label{e.1.2}
g|_{T(\partial M)} = \gamma .
\end{equation}

 Although there has been a great deal of interest in such existence (and 
uniqueness) questions on compact manifolds without boundary, very little in 
the way of general results or a general theory are known, cf.~\cite{B, LW} for 
surveys. Similarly, this question has been extensively studied for complete 
metrics on non-compact manifolds, particularly in the asymptotically Euclidean, 
flat and asymptotically hyperbolic settings. However, Einstein metrics on manifolds 
with boundary, which are in a sense intermediate between the compact and complete, 
non-compact cases, have not been studied in much detail in the literature. 

\medskip

 To describe the results, for a given $\lambda\in{\mathbb R}$, let ${\mathcal E}  = 
{\mathcal E}_{\lambda}^{m,\alpha}(M)$ be the moduli space of Einstein metrics on $M$, 
satisfying (1.1), which are $C^{m,\alpha}$ smooth up to $\partial M$; here 
$m \geq 3$ and $\alpha \in (0,1)$. By definition, ${\mathcal E}$ is the space 
of all such metrics satisfying (1.1), modulo the action of the group 
${\mathcal D}_{1} = {\mathcal D}_{1}^{m+1,\alpha}$ of $C^{m+1,\alpha}$ 
diffeomorphisms of $M$ equal to the identity on $\partial M$. 

  The first main result of the paper is the following:

\begin{theorem} \label{t1.1} 
Suppose $\pi_{1}(M, \partial M) = 0$. Then for any $\lambda \in {\mathbb R}$, the 
moduli space ${\mathcal E}$, if non-empty, is an infinite dimensional $C^{\infty}$ 
smooth Banach manifold. 
\end{theorem}

  Theorem 1.1 also holds in the $C^{\infty}$ context: the space ${\mathcal E}^{\infty}$ 
of $C^{\infty}$ Einstein metrics on $M$ is a smooth Fr\'echet manifold.

  The topological condition $\pi_{1}(M, \partial M) = 0$ means that $\partial M$ is 
connected, and the inclusion map $\iota: \partial M \rightarrow M$ induces a surjection 
$$\pi_{1}(\partial M) \rightarrow \pi_{1}(M) \rightarrow 0.$$
It is an open question whether Theorem 1.1 holds without this topological condition. 
The method of proof, via the implicit function theorem, fails without it, cf.~Remark 
2.6. On the other hand, Theorem 1.1 holds at least for generic Einstein metrics, 
without the $\pi_{1}$ condition, if $\partial M$ is connected. 

  A consequence of the proof of Theorem 1.1 is that the moduli space ${\mathcal E}$ 
is ``unobstructed'', in that any infinitesimal Einstein deformation $h$ of $(M, g)$ 
is tangent to curve in ${\mathcal E}$, i.e.~all infinitesimal deformations may be 
integrated to curves. This is in strong contrast to the situation on compact 
manifolds without boundary, where a well-known result of Koiso~\cite{K} gives examples 
where the Einstein moduli space is obstructed, cf.~also ~\cite{B}. 

\medskip

  Theorem 1.1 does not involve the specification of any boundary values of the 
metric $g$. Boundary values are given by natural boundary maps to the space 
of symmetric bilinear forms $S_{2}(\partial M)$ on $\partial M$. For the Dirichlet 
problem, one has the $C^{\infty}$ smooth Dirichlet boundary map
\begin{equation} \label{e1.3}
\Pi_{D}:{\mathcal E}^{m,\alpha} \rightarrow  Met^{m,\alpha}(\partial M), \ \ 
\Pi_{D}[g] = \gamma = g|_{T(\partial M)},
\end{equation}
where $Met^{m,\alpha}(\partial M)$ is the Banach space of $C^{m,\alpha}$ metrics on 
$\partial M$. However, $\Pi_{D}$ does not have good local properties, in that 
$\Pi_{D}$ is never Fredholm. For instance, when $m < \infty$, $D\Pi$ always has 
an infinite dimensional cokernel, so that the variety ${\mathcal B} = 
\Pi({\mathcal E}^{m,\alpha})$ has infinite codimension in $Met^{m,\alpha}(\partial M)$. 
This is a consequence of the scalar or Hamiltonian constraint on the boundary 
metric $\gamma$ induced by the Einstein metric $(M, g)$:
\begin{equation}\label{e1.4}
|A|^{2} - H^{2} + s_{\gamma} - (n-1)\lambda = 0.
\end{equation}
Here $A$ is the $2^{\rm nd}$ fundamental form of $\partial M$ in $(M, g)$, 
$H = tr A$ is the mean curvature and $s_{\gamma}$ is the scalar curvature of 
$(\partial M, \gamma)$. For $g \in {\mathbb E}^{m,\alpha}$, one has 
$A, H \in S_{2}^{m-1,\alpha}(\partial M)$, so that (1.4) gives $s_{\gamma} \in 
C^{m-1,\alpha}(\partial M)$. However, a generic $C^{m,\alpha}$ metric $\gamma$ on 
$\partial M$ has scalar curvature $s_{\gamma}$ in $C^{m-2,\alpha}$; in fact the 
space of $C^{m,\alpha}$ metrics $\gamma$ on $\partial M$ for which $s_{\gamma} 
\in C^{m-1,\alpha}(\partial M)$ is of infinite codimension. Of course the simplest 
instance of this relation is Gauss' Theorema Egregium, $K = \frac{1}{2}s_{\gamma} 
= det A$, for surfaces in ${\mathbb R}^{3}$. 

  Similarly, there are situations where the linearization $D\Pi$ has infinite 
dimensional kernel; for example this is the case whenever the $2^{\rm nd}$ 
fundamental form $A$ of $\partial M$ in $M$ vanishes on an open set in $\partial M$. 
These remarks show that the Dirichlet problem for the Einstein equations is not a 
well-posed elliptic boundary value problem. The discussion above also holds for 
the natural Neumann boundary map, taking $g \in {\mathcal E}$ to its $2^{\rm nd}$ 
fundamental form $A$ on $\partial M$. 

  These failures of the Fredholm property above are closely related to the fact that 
Einstein metrics are invariant under the full diffeomorphism group ${\mathcal D}$ 
of $\bar M$, which is much larger than the restricted group ${\mathcal D}_{1}$. 
It is also closely related to loss-of-derivative issues in the isometric 
embedding of manifolds in ${\mathbb R}^{N}$, cf.~\cite{Na, Ni}. 

   On the other hand, there are Fredholm boundary maps of mixed (Dirichlet-Neumann) 
type. There are several classes of these, but perhaps the most natural is given 
by the following result. Let ${\mathcal C}^{m,\alpha}(\partial M)$ be the space 
of pointwise conformal classes of $C^{m,\alpha}$ metrics on $\partial M$. 

\begin{theorem}\label{t1.2}
The boundary map
\begin{equation}\label{e1.5}
\widetilde \Pi_{D}: {\mathcal E}^{m,\alpha} \rightarrow 
{\mathcal C}^{m,\alpha}(\partial M)\times C^{m-1,\alpha}(\partial M),
\end{equation}
$$\widetilde \Pi_{D}(g) = ([\gamma], H),$$
is $C^{\infty}$ smooth and Fredholm, of Fredholm index 0.
\end{theorem}

   In particular, the image $\widetilde{\mathcal B} = \widetilde \Pi_{D}
({\mathcal E}^{m,\alpha})$ is a variety of finite codimension in 
${\mathcal C}^{m,\alpha}(\partial M)\times C^{m-1,\alpha}(\partial M)$. 
It is an interesting open problem to relate this image with the image of 
the usual (non-Fredholm) Dirichlet boundary map \eqref{e1.3}. Thus, one may 
fix the conformal class $[\gamma]$ and vary the mean curvature $H$. It would be 
interesting to understand the resulting space of metrics $\widetilde{\mathcal B} 
\cap [\gamma]$ within $[\gamma]$ that are obtained in this way. 

\medskip

  The results above generalize easily to ``exterior'' boundary value problems. 
In this context, $(M, g)$ is then a complete, non-compact manifold, with a compact 
(interior) boundary. Such Einstein metrics necessarily have non-positive scalar 
curvature, and the simplest asymptotic behaviors are asymptotically (locally) 
hyperbolic, when $s < 0$, and asymptotically (locally) Euclidean or flat. The 
former case has been extensively studied elsewhere, (cf.~\cite{An1} for example), 
in the case $\partial M = \emptyset$, so we concentrate here on the asymptotically 
flat case. 

 Suppose then $M$ is a manifold with a compact non-empty boundary $\partial M$ 
and a finite collection of asymptotically locally flat ends; such ends are 
metrically asymptotic to a flat metric on the space $({\mathbb R}^{m} \times 
T^{n+1-m})/ \Gamma$, where $T^{k}$ is the $k$-torus, $1 \leq m \leq n+1$ and 
$\Gamma$ is a finite group of Euclidean isometries. 

\begin{theorem}\label{t1.3}
The results above, i.e.~Theorems {\rm 1.1-1.2}, hold for the moduli space 
${\mathcal E}$ of Ricci-flat, locally asymptotically flat metrics on $M$, 
\end{theorem}

  A more precise statement of Theorem 1.3, in particular regarding the assumptions 
on the asymptotic behavior of the metrics, is given in \S 4, cf.~Theorem 4.2. 

\medskip

  The results above also hold for the Einstein equations coupled to other fields, 
for example scalar fields, sigma models (harmonic maps), etc. These are discussed 
in detail in \S 5. In fact the method of proof is quite general and should apply 
to many geometric variational problems. 

\medskip

  Theorems 1.1 and 1.2 show that one has reasonably good local behavior associated 
with the moduli space ${\mathcal E}$ of Einstein metrics on $M$. It is then of 
basic interest to understand more global issues associated with ${\mathcal E}$; 
for example, under what conditions is the boundary map $\widetilde \Pi_{D}$ in 
\eqref{e1.5} proper? We hope to address some of these questions in the future.

\section{The Moduli Space ${\mathcal E}$.}
\setcounter{equation}{0}

  Theorem 1.1 is proved via application of the implicit function theorem, 
(i.e.~the regular value theorem). To do this, one needs to choose suitable 
function spaces and make a choice of gauge in order to break the diffeomorphism 
invariance of the Einstein equations. As function space, we consider the Banach space 
\begin{equation} \label{e2.1}
Met(M) = Met^{m,\alpha}(M) 
\end{equation}
of metrics on $M$ which are $C^{m,\alpha}$ smooth up to $\partial M$. 
Here $m$ is any fixed integer with $m \geq 2$ and $\alpha \in (0,1)$. 
In the following, the smoothness index $(m, \alpha)$ will often be 
suppressed from the notation unless it is important to indicate it. Let 
\begin{equation} \label{e2.2}
{\mathbb E} = {\mathbb E}(M)
\end{equation}
be the space of Einstein metrics on $M$, 
\begin{equation} \label{e2.3}
Ric_{g} = \lambda g, 
\end{equation}
viewed as a subset of $Met(M)$, for any fixed $\lambda \in {\mathbb R}$. The 
Einstein operator $E$ is a ($C^{\infty}$) smooth map
\begin{equation}\label{e2.4}
E: Met (M) \rightarrow S_{2}(M),
\end{equation}
$$E(g) = Ric_{g} - \lambda g,$$
or more precisely $E: Met^{m,\alpha}(M) \rightarrow  
S_{2}^{m-2,\alpha}(M)$, where $S_{2}^{m-2,\alpha}(M)$ is the space of 
$C^{m-2,\alpha}$ symmetric bilinear forms on $M$. Thus
$${\mathbb E} = E^{-1}(0).$$
 
 Let $\widetilde g \in {\mathbb E}$ be a fixed but arbitrary background Einstein metric. 
A number of different gauge choices have been used to study the Einstein equations 
(2.3) near $\widetilde g$. For the purposes of this work, the simplest and most 
natural choice is the Bianchi-gauged Einstein operator, given by 
\begin{equation} \label{e2.5}
\Phi_{\widetilde g}: Met(M) \rightarrow  S_{2}(M), 
\end{equation}
$$\Phi_{\widetilde g}(g) = Ric_{g} - \lambda g + 
\delta_{g}^{*}\beta_{\widetilde g}(g),$$
where $(\delta^{*}X)(A,B) = \frac{1}{2}(\langle \nabla_{A}X, B\rangle  + 
\langle \nabla_{B}X, A\rangle )$ and $\delta X = -tr \delta^{*}X$ is the 
divergence and $\beta_{\widetilde g}(g) = \delta_{\widetilde g} g + 
\frac{1}{2}d tr_{\widetilde g} g$ is the Bianchi operator with respect to 
$\widetilde g$. Although $\Phi_{\widetilde g}$ is defined for all $g\in Met(M)$, 
we will only consider it acting on $g$ near $\widetilde g$. 

 Clearly $g$ is Einstein if $\Phi_{\widetilde g}(g) = 0$ and 
$\beta_{\widetilde g}(g) = 0$, so that $g$ is in the Bianchi-free 
gauge with respect to $\widetilde g$. Using standard formulas for the linearization 
of the Ricci and scalar curvatures, cf.~\cite{B} for instance, one finds that the 
linearization of $\Phi$ at $\widetilde g = g$ is given by 
\begin{equation} \label{e2.6}
L(h) = 2(D\Phi_{\widetilde g})_{g}(h) = D^{*}Dh - 2Rh 
\end{equation}
where the covariant derivatives and curvature are taken with respect to 
$g$. Similarly, the linearization $E' = L_{E}$ of the Einstein operator $E$ is 
given by 
\begin{equation} \label{e2.7}
2L_{E}(h) = L(h) - 2\delta^{*}\beta (h).
\end{equation}
Note that the operator $L$ is formally self-adjoint. While $L$ is elliptic, 
$L_{E}$ is not; this is the reason for choosing a gauge. The zero-set of 
$\Phi_{\widetilde g}$ near $\widetilde g$, 
\begin{equation} \label{e2.8}
Z = \{g: \Phi_{\widetilde g} = 0\}, 
\end{equation}
consists of metrics $g \in Met(M)$ satisfying the Ricci soliton equation
$$Ric_{g} - \lambda g + \delta_{g}^{*}\beta_{\widetilde g}(g) = 0.$$

\medskip

 One needs to choose boundary conditions on $\partial M$ to obtain a 
well-defined elliptic boundary value problem for the operator $\Phi$ 
on $M$. This will be done in detail in \S 3. For now, given $\widetilde g$, 
consider simply the Banach space
\begin{equation}\label{e2.9}
Met_{C}(M) = Met_{C}^{m,\alpha}(M) = \{g \in Met^{m,\alpha}(M): 
\beta_{\widetilde g}(g) = 0 \ {\rm on} \ \partial M \}. 
\end{equation}
We only consider metrics $g \in Met_{C}(M)$ near the background $\widetilde g$. 
Clearly the map 
$$\Phi : Met_{C}(M) \rightarrow  S^{2}(M),$$
is $C^{\infty}$ smooth. 

  Let $Z_{C}$ be the space of metrics $g\in Met_{C}(M)$ satisfying 
$\Phi_{\widetilde g}(g) = 0$, and let 
\begin{equation} \label{e2.10}
{\mathbb E}_{C} \subset  Z_{C} 
\end{equation}
be the subset of Einstein metrics $g$, $Ric_{g} = \lambda g$ in 
$Z_{C}$. Next we need to show that the opposite inclusion to \eqref{e2.10} 
holds, so that ${\mathbb E}_{C} = Z_{C}$. Let $\chi_{1}^{k,\alpha}$ be the 
space of $C^{k,\alpha}$ vector fields on $M$ which vanish on $\partial M$. 
One then has
$$V = \beta_{\widetilde g}(g) \in \chi_{1}^{m-1,\alpha},$$
and one needs to show that $\delta^{*}V = 0$. (Here and below we identify 
vector fields and 1-forms via the metric $g$). This will require several 
Lemmas, which will also be of importance later. 
\begin{lemma}\label{l2.1}
For $g$ in $Met^{m,\alpha}(M)$, one has
\begin{equation}\label{e2.11}
T_{g}Met^{m-2,\alpha}(M) \simeq S_{2}^{m-2,\alpha}(M) = Ker \delta \oplus 
Im \delta^{*},
\end{equation}
where $\delta^{*}$ acts on $\chi_{1}^{m-1,\alpha}$. 
\end{lemma}

{\bf Proof:} Given $h \in S_{2}^{m-2,\alpha}(M)$, consider the equation 
$\delta \delta^{*}X = \delta h \in C^{m-3,\alpha}$. If $X = 0$ at 
$\partial M$, this has a unique solution $X$ with $X \in \chi_{1}^{m-1,\alpha}$, 
by elliptic regularity ~\cite{GT, Mo}. Setting $\pi = h - \delta^{*}X$ gives the 
splitting \eqref{e2.11}. 
{\endproof}

\begin{lemma}\label{l2.2}
For $\widetilde g \in {\mathbb E}^{m,\alpha}$ and $g$ in $Met^{m,\alpha}$ close 
to $\widetilde g$, one has
\begin{equation}\label{e2.12}
T_{g}Met^{m-2,\alpha}(M) \simeq S_{2}^{m-2,\alpha}(M) = Ker \beta \oplus 
Im \delta^{*}.
\end{equation}
\end{lemma}

{\bf Proof:} By the same argument as in Lemma 2.1, it suffices to prove that 
the operator $\beta \delta^{*}: \chi_{1}^{m-1,\alpha} \rightarrow \Omega^{1}$ 
is an isomorphism, where $\Omega^{1}$ is the space of $C^{m-3,\alpha}$ 1-forms 
on $M$. Since this is an open condition, it suffices to prove this when 
$g = \widetilde g$ is Einstein. A standard Weitzenbock formula gives
$$2\beta \delta^{*}X = D^{*}DX - Ric (X) =  D^{*}DX - \lambda X.$$
Hence, if $\lambda \leq 0$, $\beta \delta^{*}$ is a positive operator and 
it follows easily that $\beta \delta^{*}$ is an isomorphism, (as in the proof 
of Lemma 2.1). 

  When $\lambda > 0$, this requires some further work. First, note that $\beta$ 
itself is surjective. To see this, suppose $Y$ is a 1-form (or vector field) 
orthogonal to $Im \beta$. Then
\begin{equation}\label{e2.13}
0 = \int_{M}\langle \beta (h), Y \rangle = \int_{M}\langle h, \beta^{*}Y \rangle 
- \int_{\partial M}[h(N, Y) - {\tfrac{1}{2}}tr h \langle Y, N \rangle].
\end{equation}
Since $h$ is arbitrary, this implies $\beta^{*}Y = \delta^{*}Y + \frac{1}{2}\delta 
Y g = 0$, and hence $\delta^{*}Y = 0$. The boundary term also vanishes, which implies 
$Y = 0$ at $\partial M$. Thus, $Y$ is a Killing field vanishing on $\partial M$, and 
hence $Y = 0$, which proves the claim. 

  To prove that $\beta \delta^{*}$ is surjective, it then suffices to show that for 
any $h\in S_{2}^{m-2,\alpha}(M)$, there exists $X$ such that $\beta \delta^{*}(X) = 
\beta(h)$. Via \eqref{e2.11}, write $h = k + \delta^{*}Y$ with $\delta k = 0$. 
Then $\beta (h) = \frac{1}{2}dtr k + \beta \delta^{*}(Y)$. This shows that it 
suffices to prove $\beta \delta^{*}$ is surjective onto exact 1-forms $df$. 

   Thus, suppose there exists $f$ such that $df \perp Im \beta \delta^{*} 
= Im(D^{*}D - \lambda I)$. Arguing just as in \eqref{e2.13}, it follows that 
$d\Delta f - \lambda(df) = 0$ on $M$, with boundary condition $df = 0$ at 
$\partial M$. Hence, $\Delta f + \lambda f = const$, with $f = const$ and 
$N(f) = 0$ at $\partial M$. It then follows from unique continuation for 
Laplace-type operators that $f = const$ on $M$, and hence $\beta \delta^{*}$ 
is surjective.

  To see that $\beta \delta^{*} = D^{*}D - \lambda I$ is injective, the family 
 $D^{*}D - t\lambda I$ for $t \in [0,1]$, with boundary condition $X = 0$ on 
$\partial M$, is a curve of elliptic boundary value problems. Since the index 
is 0 when $t = 0$, it follows that the index is also 0 when $t = 1$, 
i.e.~$\beta \delta^{*}$ has index 0 on $\chi_{1}$, which proves the injectivity. 
This completes the proof. 
{\endproof}

\begin{corollary}\label{c2.3}
Any metric $g\in Z_{C}$ near $\widetilde g$ is necessarily Einstein, with 
$Ric_{g} = \lambda g$, and in Bianchi gauge with respect to $\widetilde g$, 
i.e.
\begin{equation}\label{e2.14}
\beta_{\widetilde g}(g) = 0.
\end{equation}
\end{corollary}

{\bf  Proof:} Since $g\in Z_{C}$, one has $\Phi (g) = 0$, i.e.
$$Ric_{g} - \lambda g + \delta_{g}^{*}\beta_{\widetilde g}(g) = 0.$$
The Bianchi identity $\beta_{g}(Ric_{g}) = 0$ implies
$$\beta_{g}(\delta_{g}^{*}(V)) = 0,$$
where $V = \beta_{\widetilde g}(g)$. By the constraint \eqref{e2.9}, the 
vector field $V$ vanishes on $\partial M$, so that $V \in 
\chi_{1}^{m-1,\alpha}$. It then follows from Lemma 2.2 that 
\begin{equation}\label{e2.15}
\delta^{*}V = 0,
\end{equation}
so that $g$ is Einstein. To prove the second statement, \eqref{e2.15} 
implies that $V$ is a Killing field on $(M, g)$ with $V = 0$ at $\partial M$ 
by \eqref{e2.9}. It is then standard that $V = 0$ on $M$ so that \eqref{e2.14} 
holds. 
{\endproof}

 By linearizing, the same proof shows that the infinitesimal version of 
Corollary 2.3 holds. Thus, if $k$ is an infinitesimal deformation of 
$g\in Z_{C}$, i.e.~$k \in Ker D\Phi$ and if $\beta_{\widetilde g}(g) 
= 0$, (for example $\widetilde g = g$), then $k$ is an infinitesimal 
Einstein deformation, i.e.~the variation of $g$ in the direction $k$ 
preserves \eqref{e2.3} to $1^{\rm st}$ order and $\beta (k) = 0$. The 
proof is left to the reader.

\medskip

  As mentioned above, Theorem 1.1 is proved via the implicit function theorem in 
Banach spaces. To set the stage for this, the natural or geometric Cauchy data for 
the Einstein equations \eqref{e2.3} on $M$ at $\partial M$ consist of the pair 
$(\gamma, A)$. If $k$ is an infinitesimal Einstein deformation of $(M, g)$, so 
that $L_{E}(k) = 0$, then the induced variation of the Cauchy data on 
$\partial M$ is given by 
$$k^{T} \ \ {\rm and} \ \ (A_{k}')^{T},$$
where $A_{k}' = \frac{1}{2}({\mathcal L}_{N}g)'$ is the variation of $A$ in 
the direction $k$, given by 
\begin{equation}\label{e2.16}
2A_{k}' = \nabla_{N}k + 2A\circ k - 2\delta^{*}(k(N)^{T}) - \delta^{*}(k_{00}N),
\end{equation}
where we have used the formula ${\mathcal L}_{N}k = \nabla_{N}k + 2A\circ k$. 

  It is proved in ~\cite{AH} that an Einstein metric $g$ is uniquely determined in 
a neighborhood of $\partial M$, up to isometry, by the Cauchy data $(\gamma, A)$. 
This also holds, with the same proof, for the linearized Einstein equations and 
this linearized unique continuation result will be needed in the proof of 
Theorem 1.1. 
\begin{proposition}\label{p2.4} ~\cite{AH} Given any $(M, g) \in {\mathbb E} = 
{\mathbb E}^{m,\alpha}$, $m \geq 3$, let $k$ be any infinitesimal Einstein 
deformation of $g$ such that 
\begin{equation}\label{e2.17}
k^{T} = 0 \ \ {\rm and} \ \ (A_{k}')^{T} = 0,
\end{equation}
at $\partial M$. Then there exists a $C^{m+1,\alpha}$ vector field $Z$, defined 
in a neighborhood of $V$ of $\partial M$, with $Z = 0$ on $\partial M$, 
such that on $V$,
\begin{equation}\label{e2.18}
k = \delta^{*}Z.
\end{equation}
\end{proposition}

We note that the boundary conditions \eqref{e2.17} are invariant under infinitesimal 
gauge transformations $k \rightarrow k + \delta^{*}Z$, with $Z = 0$ on $\partial M$. 

\begin{proposition}\label{p2.5}
Suppose $\pi_{1}(M, \partial M) = 0$. Then at any $\widetilde g \in {\mathbb E}^{m,\alpha}$, 
$m \geq 3$, the map $\Phi = \Phi_{\widetilde g}$ is a submersion on $Met_{C}^{m,\alpha}(M)$. 
Thus, the linearized operator $L = 2D\Phi$:
\begin{equation} \label{e2.19}
L: T_{\widetilde g}Met_{C}(M) \rightarrow  S_{2}(M) 
\end{equation}
is surjective, and the kernel of $L$ splits in $T_{\widetilde g}Met_{C}(M)$. 
\end{proposition}

{\bf Proof:} The operator $L$ is elliptic on $T_{\widetilde g}Met_{C}(M)$, and 
so Fredholm. (More precisely, one can augment the constraint \eqref{e2.9} with 
further boundary conditions to obtain an elliptic boundary value problem; this is 
discussed in detail in \S 3). 

  In particular, $Im(L)$ is closed and has a closed complement in 
$S_{2}(M)$. If $L$ is not surjective, then there exists a non-zero $k\in 
T_{\widetilde g}Met(M) = S_{2}(M)$ such that, for all $h \in T_{\widetilde g}
Met_{C}(M)$, 
\begin{equation}\label{e2.20}
\int_{M}\langle L(h), k \rangle dV_{\widetilde g} = 0.
\end{equation}
The idea of the proof is to show that \eqref{e2.20} implies $k$ is an infinitesimal Einstein 
deformation, so $L_{E}(k) = 0$, (in transverse-traceless gauge), and satisfying the 
boundary conditions \eqref{e2.17}, so that the unique continuation property in Proposition 2.4 
applies. Once this is established, the proof follows by a simple global argument, using the 
condition on $\pi_{1}$. In the following, we set $\widetilde g = g$ and drop the volume 
forms from the notation. 

  To begin, integrating \eqref{e2.20} by parts, one obtains
\begin{equation} \label{e2.21}
0 = \int_{M}\langle L(h), k \rangle = \int_{M}\langle h, L(k)\rangle + 
\int_{\partial M}D(h,k), 
\end{equation}
where the boundary pairing $D(h, k)$ has the form
\begin{equation}\label{e2.22}
D(h,k) = \langle h, \nabla_{N}k \rangle - \langle k, \nabla_{N}h \rangle.
\end{equation}

  Since $h$ is arbitrary in the interior, the bulk integral and the boundary integral 
on the right in \eqref{e2.21} vanish separately, and hence 
\begin{equation}\label{e2.23}
L(k) = 0.
\end{equation}
Next, observe that \eqref{e2.20} implies that 
\begin{equation}\label{e2.24}
\delta k = 0 \ \ {\rm on} \ \  M,
\end{equation}
so that $k$ is in divergence-free gauge on $M$. To see this, \eqref{e2.7} implies that 
$L(\delta^{*}X) = \delta^{*}Y,$ where $Y = 2\beta\delta^{*}X$. Then
$$0 = \int_{M}\langle L(\delta^{*}X), k\rangle = \int_{M}\langle Y, \delta k \rangle 
+ \int_{\partial M}k(Y, N). $$
For $h = \delta^{*}X$, the Bianchi constraint \eqref{e2.9} gives exactly $Y = 0$ at 
$\partial M$. By Lemma 2.2, $Y$ is arbitrary in the interior of $M$, which thus gives 
\eqref{e2.24}. 

  The boundary integral in \eqref{e2.21} vanishes for all $h \in T_{g}Met_{C}(M)$, 
i.e.~all $h$ satisfying the linearized constraint \eqref{e2.9}. Written out in 
tangential and normal components, this requires
\begin{equation}\label{e2.25}
(\nabla_{N}h)(N)^{T} = \delta^{T}h^{T} - \alpha(h(N)) + {\tfrac{1}{2}}d^{T}tr h,
\end{equation}
\begin{equation}\label{e2.26}
N(h_{00}) = \delta^{T}(h(N)^{T}) - h_{00}H + \langle A, h \rangle + {\tfrac{1}{2}}N(tr h),
\end{equation}
where $\alpha(h(N)) = A(h(N)) + Hh(N)^{T}$, and $N(h_{00}) = (\nabla_{N}h)(N,N)$. 

  We first use various test-forms $h$ to obtain restrictions on $k$ at $\partial M$. 
Thus, suppose first $h = 0$ at $\partial M$. The constraints \eqref{e2.25}-\eqref{e2.26} 
then require that 
$$(\nabla_{N}h)(N)^{T} = 0 \ {\rm and} \  N(h_{00}) = \langle \nabla_{N}h, \gamma 
\rangle ,$$
where we have used the fact that $tr h = h_{00} + tr_{\gamma}h$. For all such 
$h$, \eqref{e2.21}-\eqref{e2.22} gives 
$$\int_{\partial M} \langle \nabla_{N}h, k \rangle = 0,$$
and hence
$$\int_{\partial M} N(h_{00})k_{00} + {\tfrac{1}{n}}\langle 
\nabla_{N}h, \gamma \rangle \langle \gamma, k \rangle + \langle (\nabla_{N}h)_{0}^{T}, 
k_{0}^{T} \rangle = 0,$$
where $k_{0}^{T}$ is the trace-free part of $k^{T}$. This implies that
\begin{equation}\label{e2.27}
k^{T} = \phi \gamma, \ \ {\rm and} \ \ k_{00} = -{\tfrac{1}{n}}tr_{\gamma}k = -\phi,
\end{equation}
for some function $\phi$. 

  Next, set $h^{T} = h_{00} = 0$ with $h(N)^{T}$ chosen arbitrarily, and similarly 
$\nabla_{N}h = 0$ except for the two relations $(\nabla_{N}h)(N)^{T} = -\alpha(h(N)^{T})$ 
and $N(tr h) = -2\delta^{T}(h(N)^{T})$. The constraints are then satisfied and via 
\eqref{e2.27}, one has
$$\int_{\partial M}\langle \nabla_{N}h, k \rangle = -2\int_{\partial M}
\langle \alpha(h(N)), k(N)^{T} \rangle - 2{\tfrac{n-1}{n+1}}\phi \delta^{T}(h(N)^{T}),$$
while 
$$\int_{\partial M}\langle \nabla_{N}k, h \rangle = 2\int_{\partial M}\langle 
h(N)^{T}, (\nabla_{N}k)(N)^{T} \rangle.$$
Now by \eqref{e2.24} one has, (as in \eqref{e2.25}), $(\nabla_{N}k)(N, T) = 
-(\nabla_{e_{i}}k)(e_{i},T) = -e_{i}(k(e_{i},T)) + k(\nabla_{e_{i}}e_{i}, T) 
+ k(e_{i}, \nabla_{e_{i}}T) = - \langle d\phi, T \rangle - \langle \alpha(k(N)), 
T \rangle$, so that
\begin{equation}\label{e2.28}
(\nabla_{N}k)(N)^{T} = - d^{T}\phi - \alpha(k(N)).
\end{equation}
Note that $\alpha$ is symmetric: $\langle \alpha (h(N)), k(N) \rangle = \langle 
\alpha(k(N)), h(N) \rangle$. It then follows from the two equations above and the 
divergence theorem that
\begin{equation}\label{e2.29}
\phi = const.
\end{equation}
We note that, analogous to \eqref{e2.26}, \eqref{e2.24} together with \eqref{e2.27} gives
\begin{equation}\label{e2.30}
N(k_{00}) = \delta^{T}(k(N)^{T}) - k_{00}H + \langle k, A \rangle = \delta^{T}(k(N)^{T}) 
+ 2H\phi.
\end{equation}

  Next, suppose $h = 0$ except for $h_{00}$, which is chosen arbitrarily, and similarly 
$\nabla_{N}h = 0$ except for the component $\nabla_{N}h(N)^{T}$. Then \eqref{e2.25} and 
\eqref{e2.26} require
$$N(h_{00}) = -Hh_{00},$$
$$(\nabla_{N}h)(N)^{T} = {\tfrac{1}{2}}d^{T}h_{00}.$$
This gives
$$\int_{\partial M}\langle \nabla_{N}h, k \rangle = 2\int_{\partial M}\langle 
(\nabla_{N}h)^{T}, k(N)^{T} \rangle - h_{00}k_{00}H$$
$$= \int_{\partial M}\langle d^{T}h_{00}, k(N)^{T} \rangle - h_{00}k_{00}H = 
\int_{\partial M}h_{00}\delta^{T}(k(N)^{T}) - h_{00}k_{00}H,$$
while
$$\int_{\partial M}\langle h, \nabla_{N}k \rangle = \int_{\partial M}h_{00}N(k_{00}).$$
Hence
\begin{equation}\label{e2.31}
N(k_{00}) = \delta^{T}(k(N)^{T}) - k_{00}H.
\end{equation}
Via \eqref{e2.30} and \eqref{e2.27}, this implies that 
\begin{equation}\label{e2.32}
H\phi = 0,
\end{equation}
so that $\phi \equiv 0$ unless $H \equiv 0$. 

  Finally, set $h^{T} = f\gamma$, $h_{00} = -\frac{n-2}{2}f$ with the rest of $h$ set to 
0. Then setting $N(h_{00}) = nfH$ with the rest of $\nabla_{N}h$ set to 0 solves the 
constraints \eqref{e2.25}-\eqref{e2.26}, and \eqref{e2.21}-\eqref{e2.22} together with 
\eqref{e2.27} and \eqref{e2.32} then gives, 
since $f$ is arbitrary,
$$\langle \nabla_{N}k, \gamma \rangle = {\tfrac{n-2}{2}}N(k_{00}),$$
or equivalently
\begin{equation}\label{e2.33}
N(tr k) = {\tfrac{n}{2}}N(k_{00}).
\end{equation}

  We now analyse the boundary term in \eqref{e2.21} in general, using the information 
obtained above. Thus, expand the inner products in (2.33) into tangential, mixed 
and normal components. Using \eqref{e2.27}, one has $\langle\nabla_{N}h, k\rangle 
= 2\langle (\nabla_{N}h)(N)^{T}, k(N)^{T}\rangle$, since the 00 and trace components 
cancel by \eqref{e2.26} and \eqref{e2.29}. Using the constraint \eqref{e2.25}, together 
with the fact that, modulo divergence terms, $\langle \delta^{T}h^{T}, 
k(N)^{T}\rangle = \langle h^{T}, (\delta^{T})^{*}k(N)^{T}\rangle = \langle h^{T}, 
\delta^{*}(k(N)^{T})\rangle$, one has 
$$\int_{\partial M}\langle\nabla_{N}h, k\rangle = 
\int_{\partial M} 2\langle h^{T}, \delta^{*}(k(N)^{T})\rangle + 
\delta^{T}(k(N)^{T})tr h - 2\langle \alpha(h(N)), k(N) \rangle$$
$$= \int_{\partial M} 2\langle h^{T}, \delta^{*}(k(N)^{T})\rangle + 
\delta^{T}(k(N)^{T})\langle h^{T}, \gamma \rangle + \delta^{T}(k(N)^{T})h_{00}
- 2\langle \alpha(h(N)), k(N) \rangle.$$
 On the other hand, $\langle \nabla_{N}k, h \rangle = \langle (\nabla_{N}k)^{T}, 
h^{T} \rangle + 2\langle \nabla_{N}k(N)^{T}, h(N)^{T} \rangle + N(k_{00})h_{00}$. 
The middle term is computed in \eqref{e2.28} and using \eqref{e2.29}, one has 
$$\int \langle\nabla_{N}k, h\rangle = \langle (\nabla_{N}k)^{T}, h^{T} \rangle 
+ N(k_{00})h_{00} - 2\langle \alpha(k(N)), h(N)\rangle.$$
Now take the difference of these terms. Recall that $\alpha$ is symmetric and from 
\eqref{e2.16}, $(\nabla_{N}k)^{T} - 2\delta^{*}k(N)^{T} = 2(A_{k}')^{T} - 2A\circ k + k_{00}A$. 
Using also \eqref{e2.31} and \eqref{e2.32} then gives
\begin{equation}\label{e2.34}
\int_{\partial M}\langle 2(A_{k}')^{T} - \delta^{T}(k(N)^{T})\gamma  - 3\phi A, 
h^{T} \rangle = 0.
\end{equation}
Since $h^{T}$ may be chosen arbitrarily consistent with the constraints, it follows 
that the integrand is 0. Taking the $\gamma$-trace, using \eqref{e2.32}, this gives
$$N(tr k) - N(k_{00}) - (n-2)\delta^{T}(k(N)^{T}) = 0,$$
which via \eqref{e2.31} and \eqref{e2.33} implies that
\begin{equation}\label{e2.35}
N(tr k) = N(k_{00}) = 0. 
\end{equation}
Finally, we claim that 
\begin{equation}\label{e2.36}
tr k = (n-1)\phi = 0. 
\end{equation}
To see this, the trace of \eqref{e2.23} gives
$$\Delta trk + \frac{s}{n+1}trk = 0.$$
Since $trk = (n-1)\phi = const$ and $N(tr k) = 0$ on $\partial M$, a standard 
unique continuation principle for the Laplacian implies that $tr k = c$ on $M$. 
However, integrating the equation above over $M$ and using \eqref{e2.35} implies 
$tr k$ has mean value 0 on $M$, which gives \eqref{e2.36}. 

  The results above thus imply that 
\begin{equation}\label{e2.37}
k^{T} = 0, \ \ (A_{k}')^{T} = 0, \ {\rm on} \ \partial M.
\end{equation}
In addition, $k$ is an infinitesimal Einstein deformation, since $\beta (k) = 0$ 
and \eqref{e2.23} holds. By the unique continuation property, one thus has
\begin{equation}\label{e2.38}
k = \delta^{*}Z,
\end{equation}
in a neighborhood $V$ of $\partial M$, with $Z = 0$ at $\partial M$. It follows 
from the topological condition $\pi_{1}(M, \partial M) = 0$ and analytic 
continuation in the interior that \eqref{e2.38} holds globally on $M$. Since $k$ 
is divergence-free by \eqref{e2.24}, one has
$$\delta \delta^{*}Z = 0,$$
globally on $M$, with $Z = 0$ on $\partial M$. Pairing this with $Z$ and integrating 
over $M$, it follows from the divergence theorem that $\delta^{*}Z = 0$ on $M$, and 
hence $k = 0$, which completes the proof of surjectivity. 

\medskip

 It is now essentially standard or formal that the kernel of 
$D\Phi_{g}$ splits, i.e.~it admits a closed complement in 
$T_{g}Met_{C}$. In more detail, it suffices to find a bounded linear 
projection $P$ mapping $T_{g}Met_{C}(M)$ onto $Ker(D\Phi_{g})$. To do 
this, let $Met_{C}^{0}(M) \subset  Met_{C}(M)$ be the subspace of 
metrics $g$ such that $g|_{T(\partial M)} = \widetilde g|_{T(\partial M)} 
= \widetilde \gamma$. Choose a fixed smooth extension operator 
taking metrics $\gamma$ on $T(\partial M)$ into $Met_{C}(M)$, and let 
$Met_{C}^{1}(M)$ be the resulting space of metrics, so that 
$$Met_{C}(M) = Met_{C}^{0}(M) \oplus  Met_{C}^{1}(M).$$
Let $T$ and $T^{0}, T^{1}$ denote the corresponding tangent spaces at 
$\widetilde g$ and $L^{i} = D\Phi|_{T^{i}}: T^{i} \rightarrow S_{2}(M)$, 
for $i = 1,2$. Then
\begin{equation} \label{e2.39}
KerD\Phi_{g} = \{(h, g_{\dot \gamma})\in T^{0}\oplus T^{1}: 
L^{0}(h) + L^{1}(g_{\dot \gamma}) = 0\}.  
\end{equation}
The operator $L$ is elliptic, so that $L^{0}$ is Fredholm on $T^{0}$. 
The image $Im(L^{0})$ thus has a finite dimensional complement $S$, 
$S_{2}(M) = Im(L^{0})\oplus S$. By \eqref{e2.39}, $Im(L^{1}) \subset Im(L^{0})$ 
and so $Im L^{1} \subset  Ker(\pi_{S}L^{1})$, where $\pi_{S}$ is orthogonal 
projection onto $S$. By the nondegeneracy property \eqref{e2.20}, $L^{1}$ maps 
onto $S$ and hence $Im \pi_{S}L^{1} = S$. Viewing $S$ as a subspace of $T$, under 
the natural isomorphism $T \simeq S_{2}(M)$, this gives $T = Im(\pi_{S}L^{1}) 
\oplus  Ker(\pi_{S}L^{1})$, i.e.~$Ker(\pi_{S}L^{1})$ splits, and so there is 
a bounded linear projection $P_{1}$ onto $Ker(\pi_{S}L^{1})$. The mapping 
$L + \pi_{S}$ is invertible and 
$$P(h, g_{\dot \gamma}) = ((L^{0}+\pi_{S})^{-1}(- L^{1}P_{1}(g_{\dot \gamma}) + 
\pi_{S}h),P_{1}g_{\dot \gamma})$$
gives required bounded linear projection onto $Ker D\Phi_{g}$. This completes 
the proof. 

{\endproof}

\begin{remark}\label{r2.6}
 {\rm There exist at least some examples of Einstein metrics $(M, g)$ having 
non-zero solutions of \eqref{e2.20}, so that $L$ on $T_{g}Met_{C}(M)$ is not 
surjective in general; (of course such metrics must violate the condition $\pi_{1}(M, 
\partial M) = 0$). As a simple example, let $M = I\times T^{n}$, $I = [0,1]$, 
and let $g$ be a flat product metric on $M$. Let $x^{\alpha}$ denote standard 
coordinates on $M$, with $x^{0} = t$ parametrizing $I$. Then the symmetric form 
\begin{equation} \label{e2.40}
k = dt\cdot  dx^{\alpha}  = \delta^{*}(t\nabla x^{\alpha} ), 
\end{equation}
is a divergence-free deformation of the flat metric satisfying \eqref{e2.20}, 
at least when $\alpha > 0$. Note that this solution is pure gauge, 
$k = \delta^{*}Y$, with $Y = t\nabla x^{\alpha}$ vanishing at one boundary 
component but not at the other. 

   On the other hand, the condition on $\pi_{1}$ in Proposition 2.5 and 
Theorem 1.1 is used only to extend the locally defined solution $Z$ in 
\eqref{e2.38} to a globally defined vector field on $M$ with $Z = 0$ on 
$\partial M$. For example, $Z$ in \eqref{e2.38} is unique modulo local Killing 
fields. Hence if $(M, g)$ has no local Killing fields, (which is the case 
for generic metrics), then Proposition 2.5, and Theorem 1.1, hold near 
$g$, provided $\partial M$ is connected. }
\end{remark}

\begin{corollary}\label{c2.7} 
Under the assumptions of Proposition 2.5, if ${\mathbb E}$ is non-empty, then 
the local spaces ${\mathbb E}_{C}$ are infinite dimensional $C^{\infty}$ Banach 
manifolds, with 
\begin{equation} \label{e2.41}
T_{\widetilde g}{\mathbb E}_{C} = Ker(D\Phi_{\widetilde g})_{\widetilde g}.
\end{equation}
\end{corollary}

{\bf Proof:}
 This is an immediate consequence Corollary 2.3, Proposition 2.5 and the 
implicit function theorem, (regular value theorem), in Banach spaces. 
{\endproof}

   By Corollary 2.4, Einstein metrics in ${\mathbb E}_{C}$ satisfy the Bianchi 
gauge condition
\begin{equation}\label{e2.42}
\beta_{\widetilde g}(g) = 0. 
\end{equation}

  We need to show that \eqref{e2.42} is actually a well-defined gauge condition. 
Let ${\mathcal D}_{1} = {\mathcal D}^{m+1,\alpha}_{1}(M)$ be the group of 
$C^{m+1,\alpha}$ diffeomorphisms of $M$ which equal the identity on $\partial M$. 
The action of ${\mathcal D}_{1}$ on ${\mathbb E}$ is continuous and also free, 
since any isometry $\phi$ of a metric inducing the identity on $\partial M$ 
must itself be the identity. However, the action of ${\mathcal D}_{1}$ is not 
apriori smooth. Namely, as before let $\chi_{1} = \chi_{1}^{m+1,\alpha}$ denote 
the space of $C^{m+1,\alpha}$ vector fields $X$ on $\bar M$ with $X = 0$ on 
$\partial M$, so that $\chi_{1}$ represents the tangent space of ${\mathcal D}_{1}$ 
at the identity. For $X \in \chi_{1}$ and $g \in Met^{m,\alpha}(M)$, one has 
$\delta^{*}X = \frac{1}{2}{\mathcal L}_{X}g \in S_{2}^{m-1,\alpha}(M)$, but 
$\delta^{*}X \notin S_{2}^{m,\alpha}(M)$, so that there is a loss of one 
derivative. 

  For the same reasons, at a general metric $g \in Met^{m,\alpha}$, the splitting 
\eqref{e2.12} does not hold when $(m-2)$ is replaced by $m$. However, for Einstein metrics, 
this loss of regularity can be restored. 
\begin{lemma}\label{l2.8}
For $g \in {\mathbb E}^{m,\alpha}$, the splittings {\rm (2.11)} and 
{\rm (2.12)} hold, and for any $X \in \chi_{1}^{m+1,\alpha}$, $\delta^{*}X 
\in S_{2}^{m,\alpha}(M)$. 
\end{lemma}

{\bf Proof:} 
By the proofs of \eqref{e2.11} and \eqref{e2.12}, it suffices to prove the second 
statement. Since Einstein metrics are $C^{\infty}$ smooth, (in fact real-analytic), 
in harmonic coordinates in the interior, ${\mathcal L}_{X}g$ is $C^{m,\alpha}$ 
smooth in the interior of $M$. To see that ${\mathcal L}_{X}g$ is $C^{m,\alpha}$ 
smooth up to $\partial M$, recall that in suitable 
boundary harmonic coordinates, one has
$$\Delta_{g}g_{\alpha\beta} + Q_{\alpha\beta}(g, \partial g) = -2Ric_{\alpha\beta} = 
-2 \lambda g_{\alpha\beta},$$
cf.~\cite{AT} for the analysis of boundary regularity of Einstein metrics. Applying 
$X$ to this equation and commuting derivatives gives an equation for 
$\Delta_{g} X(g_{\alpha\beta})$ with 0 boundary values, (since $X(g) = 0$ on 
$\partial M$), and with right-hand side in $C^{m-2,\alpha}$. Elliptic boundary 
regularity results then imply that $X(g_{\alpha\beta}) \in C^{m,\alpha}$. From this, 
it is easy to see that ${\mathcal L}_{X}g$ is $C^{m,\alpha}$ smooth up to $\partial M$. 
{\endproof}

  Next we pass from the infinitesimal splitting to its local version. 

\begin{lemma}\label{l2.9} 
Given any $\widetilde g \in {\mathbb E}^{m,\alpha}$ and $g \in Met^{m,\alpha}(M)$ 
nearby $\widetilde g$, there exists a unique diffeomorphism $\phi \in 
{\mathcal D}_{1}^{m+1,\alpha}$, close to the identity, such that 
\begin{equation} \label{e2.43}
\beta_{\widetilde g}(\phi^{*}g) = 0. 
\end{equation}
In particular, $\phi^{*}g\in Met_{C}^{m,\alpha}(M)$.
\end{lemma}

{\bf Proof:} Given Lemma 2.8, this can be derived from the slice theorem of Ebin ~cite{Eb}, 
but we give a direct and simpler argument here. (Note that the Lemma does not assert 
the existence of a smooth slice). Let $\widetilde g\in {\mathbb E}$ and consider the 
map $F: {\mathcal D}_{1}\times Met_{C}(M) \rightarrow Met(M)$ given by 
$F(\phi, g) = \phi^{*}g$. The proof of Lemma 2.8 above shows that $F$ is 
$C^{\infty}$ smooth at $\widetilde g$; $F$ is linear in $g$ and smooth in the direction 
of ${\mathcal D}_{1}$ at $\widetilde g$ and hence smooth at $\widetilde g$. 

 Now suppose $g \in Met(M)$ is close to $\widetilde g$. The linearization of 
$F$ at $(Id, \widetilde g)$ is the map $(X, h) \rightarrow  \delta^{*}X 
+ h$. By Lemma 2.8, given any $h$, there exists a unique vector field 
$X \in \chi_{1}^{m+1,\alpha}$ such that on $(M, \widetilde g)$, 
$\beta (\delta^{*}X + h) = 0$, or $$\beta\delta^{*}X = -\beta (h), $$
with respect to $\widetilde g$. Hence, for $g$ sufficiently close to $\widetilde g$, 
there is a vector field $X$ such that
$$|\beta_{\widetilde g}(g + \delta_{\widetilde g}^{*}X)| << \beta_{\widetilde g}(g).$$
It then follows from the inverse function theorem that there 
exists a unique diffeomorphism $\phi\in{\mathcal D}_{1}$ close to the identity such 
that \eqref{e2.43} holds. 
{\endproof}

  Lemma 2.9 implies that if $g \in {\mathbb E}$ is an Einstein metric near 
$\widetilde g$, then $g$ is isometric, by a unique diffeomorphism in ${\mathcal D}_{1}$, 
to an Einstein metric in ${\mathbb E}_{C}$. Hence \eqref{e2.42} is a well-defined gauge 
condition and the spaces ${\mathbb E}_{C}$ are local slices for the action of 
${\mathcal D}_{1}$ on ${\mathbb E}$. 

 We are now in position to complete the proof of Theorem 1.1. 

{\bf Proof of Theorem 1.1.}

 The space ${\mathbb E} = {\mathbb E}^{m,\alpha}(M) \subset Met^{m,\alpha}(M)$ of 
all Einstein metrics on $M$ is invariant under the action of the group 
${\mathcal D}_{1} = {\mathcal D}_{1}^{m+1,\alpha}$. The moduli space ${\mathcal E}  
= {\mathcal E}^{m,\alpha}(M)$ of $C^{m,\alpha}$ Einstein metrics on $M$ is the quotient
$${\mathcal E}  = {\mathbb E}/{\mathcal D}_{1}.$$
Two metrics $g_{1}$ and $g_{2}$ in ${\mathcal E}$ are equivalent if there exists 
$\phi \in {\mathcal D}_{1}$, such that $\phi^{*}g_{1} = g_{2}$. 

  The local spaces ${\mathbb E}_{C}$ are smooth Banach manifolds and depend smoothly 
on the background metric $\widetilde g$, since the gauge condition \eqref{e2.42} varies 
smoothly with $\widetilde g$. As noted above, the action of ${\mathcal D}_{1}$ on 
${\mathbb E}$ is free and by Lemma 2.8, the action is smooth. Hence the global space 
${\mathbb E}$ is a smooth Banach manifold, as is the quotient ${\mathcal E}$. The local 
slices ${\mathbb E}_{C}$ represent local coordinate patches for ${\mathcal E}$. 

  It also follows immediately from the proof above that the spaces ${\mathbb E}^{\infty}$ 
and ${\mathcal E}^{\infty} = {\mathbb E}^{\infty} / {\mathcal D}^{\infty}$ of $C^{\infty}$ 
Einstein metrics on $M$ are smooth Fr\'echet manifolds. 
{\endproof}

\section{Elliptic Boundary Problems for the Einstein Equations.}
\setcounter{equation}{0}

  In this section, we consider elliptic boundary value problems for the 
Einstein equations. We begin with the Dirichlet boundary value problem. 
A metric $g$ on $M$ induces naturally a boundary metric
\begin{equation}\label{e3.1}
\gamma  = g^{T} = g|_{T(\partial M)}
\end{equation}
on $\partial M$. One also has a normal part $g^{N} \equiv g|_{N(\partial M)}$ 
of the metric $g$ at $\partial M$, i.e. the restriction of $g$ to the normal 
bundle of $\partial M$ in $M$. In local coordinates $(x^{0}, x^{1}, \dots , x^{n})$ 
for $\partial M$ in $M$ with $x^{0} = 0$ on $\partial M$, these are the 
$g_{0\alpha}$ components of $g_{\alpha\beta}$, with $0 \leq \alpha \leq n$. 
Observe that the normal part of $g$ is a gauge term, in the sense that it 
transforms as a 1-form under the action of diffeomorphisms of $M$ equal to 
the identity on $\partial M$.

  Given the work in \S 2 and the relation \eqref{e2.7} between $L_{E}$ with $L$, 
the most obvious boundary conditions to impose for the Dirichlet problem are:
\begin{equation} \label{e3.2}
g|_{T(\partial M)} = \gamma  \ {\rm on} \ \partial M, \ \ {\rm and} 
\end{equation}
\begin{equation} \label{e3.3}
\beta_{\widetilde g}(g) = 0 \ {\rm on} \  \partial M. 
\end{equation}
Here $\gamma$ is an arbitrary Riemannian metric on $\partial M$, close 
to $\widetilde \gamma$ in $Met^{m,\alpha}(\partial M)$. Note this is a 
formally determined set of boundary conditions; the Dirichlet condition 
\eqref{e3.2} gives $\frac{1}{2}n(n+1)$ equations, while the Neumann-type 
boundary condition \eqref{e3.3} gives $n+1$ equations. In sum, this gives 
$\frac{1}{2}(n+1)(n+2)$ equations, which equals the number of 
components of the variable $g$ on $M$. 

  However, the operator $\Phi$ with the boundary conditions \eqref{e3.2}-\eqref{e3.3} 
does not form a well-defined elliptic boundary value problem. Geometrically, the 
reason for this is as follows. Metrics $g$ satisfying $\Phi(g) = 0$ with the 
boundary condition \eqref{e3.3} are Einstein, (cf.~Corollary 2.4), and so satisfy 
the Einstein constraint equations on $\partial M$. These are given by 
\begin{equation}\label{e3.4}
|A|^{2} - H^{2} + s_{\gamma} - (n-1)\lambda = (Ric_{g} - \lambda g)(N,N) = 0,
\end{equation}
\begin{equation}\label{e3.5}
\delta(A - H\gamma) = Ric(N, \cdot) = 0.
\end{equation}
The scalar or Hamiltonian constraint \eqref{e3.5} imposes a constraint on the 
regularity of the boundary metric $\gamma$ not captured by \eqref{e3.2}-\eqref{e3.3}. 
Thus, if the boundary conditions \eqref{e3.2}-\eqref{e3.3} gave an elliptic system, 
\eqref{e3.4} would hold for a space of boundary metrics $\gamma$ of finite codimension 
in $Met^{m,\alpha}(\partial M)$, which, as discussed in the Introduction, is impossible.

  The discussion above implies there is no natural elliptic boundary value problem 
for the Einstein equations, associated with Dirichlet boundary values. To obtain 
an elliptic problem, one needs to add either gauge-dependent terms or terms 
depending on the extrinsic geometry of $\partial M$ in $(M, g)$. To maintain a 
determined boundary value problem, one then has to subtract part of the 
intrinsic Dirichlet boundary data on $\partial M$. 

  There are several ways to carry this out in practice, but we will concentrate 
on the following situations. Let $B$ be a $C^{m,\alpha}$ positive definite 
symmetric bilinear form on $\partial M$. In place of prescribing the boundary 
metric $g^{T}$ on $\partial M$, only $g^{T}$ modulo $B$ will be prescribed. Thus, 
let $\pi$ be the projection
$$\pi: Met^{m,\alpha}(\partial M) \rightarrow Met^{m,\alpha}(\partial M)/ B, 
\ \ \pi(\gamma) = [\gamma]_{B} = [\gamma + fB]_{B}.$$
We allow here $B$ to depend on $\gamma$. For instance, if $B = \gamma$, then 
$[\gamma]_{B} = [\gamma]$ is the conformal class of $\gamma$.

  The simplest gauge-dependent term one can add to \eqref{e3.3} is the equation 
$g(\widetilde N, \widetilde N) = \gamma_{00}$, where $\widetilde N$ is the unit normal 
with respect to $\widetilde g$, while the simplest extrinsic geometric scalar 
is $H$, the mean curvature of $\partial M$ in $(M, g)$. 

\begin{proposition}\label{p3.1}
The Bianchi-gauged Einstein operator $\Phi$ with boundary conditions either 
\begin{equation}\label{e3.6}
\beta_{\widetilde g}(g) = 0, \ \ [g^{T}]_{B} = [\gamma]_{B}, \ \ 
g(\widetilde N, \widetilde N) = \gamma_{00} \ \ {\rm at} \ \ \partial M,
\end{equation}
or 
\begin{equation}\label{e3.7}
\beta_{\widetilde g}(g) = 0, \ \ [g^{T}]_{B} = [\gamma]_{B}, \ \ 
H_{g} = h \ \ {\rm at} \ \ \partial M,
\end{equation}
is an elliptic boundary value problem of Fredholm index 0.
\end{proposition} 

{\bf Proof:} The proof is essentially a standard computation, following ideas 
initially introduced by Nash ~\cite{Na} in the isometric embedding problem, 
cf.~also ~\cite{H}. We will follow the method used by Schlenker in ~\cite{S}. 

   It suffices to show that the leading order part of the linearized operators 
forms an elliptic system. The leading order symbol of $L = D\Phi$ is given by 
\begin{equation}\label{e3.8}
\sigma(L) = -|\xi|^{2}I,
\end{equation}
where $I$ is the $N\times N$ identity matrix, with $N = (n+2)(n+1)/2$ the dimension 
of the space of symmetric bilinear forms on ${\mathbb R}^{n+1}$. In the following, 
the subscript 0 represents the direction normal to $\partial M$ in $M$, and Latin 
indices run from $1$ to $n$. The positive roots of \eqref{e3.8} are $i|\xi|$, with 
multiplicity $N$. 

  Writing $\xi = (z, \xi_{i})$, the symbols of the leading order terms in the boundary 
operators are given by:
$$-2izh_{0k} - 2i\sum \xi_{j}h_{jk} + i\xi_{k}tr h = 0,$$
$$-2izh_{00} - 2i\sum \xi_{k}h_{0k} + iztr h = 0,$$
$$h^{T} = (\gamma')^{T} \ \ mod\, B,$$
$$h_{00} = \omega \ \ {\rm or} \ \ H_{h}' = \omega,$$
where $h$ is an $N\times N$ matrix. Then ellipticity requires that the operator defined 
by the boundary symbols above has trivial kernel when $z$ is set to the root $i|\xi|$. 
Carrying this out then gives the system 
\begin{equation}\label{e3.9}
2|\xi|h_{0k} - 2i\sum \xi_{j}h_{jk} + i\xi_{k}tr h = 0,
\end{equation}
\begin{equation}\label{e3.10}
2|\xi|h_{00} - 2i\sum \xi_{k}h_{0k} - |\xi|tr h = 0,
\end{equation}
\begin{equation}\label{e3.11}
h_{kl} = \phi b_{kk}\delta_{kl},
\end{equation}
\begin{equation}\label{e3.12}
h_{00} = 0 \ \ {\rm or} \ \ H_{h}' = 0,
\end{equation}
where without loss of generality we assume $B$ is diagonal, with entries $b_{kk}$, and 
$\phi$ is an undetermined function. 

  Multiplying \eqref{e3.9} by $i\xi_{k}$ and summing gives
$$2|\xi|i\sum \xi_{k}h_{0k} =  2i^{2}\xi_{k}^{2}h_{kk} - i^{2}\xi_{k}^{2}tr h.$$
Substituting \eqref{e3.10} on the term on the left above then gives
$$2|\xi|^{2}h_{00} - |\xi|^{2}tr h =  -2\sum \xi_{k}^{2}h_{kk} + |\xi|^{2}tr h,$$
so that
$$|\xi|^{2}h_{00} - |\xi|^{2}tr h =  -\sum \xi_{k}^{2}h_{kk} = 
-\phi \langle B(\xi), \xi \rangle.$$
Using the fact that $\sum h_{kk} = tr h - h_{00}$, this is equivalent to 
$$\phi\langle B(\xi), \xi \rangle = \phi |\xi|^{2}tr B.$$
Since $B$ is assumed to be positive definite, it follows that $\phi = 0$ and hence 
$h^{T} = 0$. 

  If the first boundary condtion $h_{00} = 0$ in \eqref{e3.12} is used, then $tr h = 0$, 
and hence via \eqref{e3.9}, $h_{0k} = 0$. This gives $h = 0$, as required. If instead 
one uses the second condition $H' = 0$ in \eqref{e3.12}, a simple computation 
shows that to leading order, $H_{h}' = tr^{T}(\nabla_{N}h - 2\delta^{*}(h(N)^{T}))$, 
which has symbol $iz\sum h_{kk} - 2i\xi_{k}h_{0k}$. Setting this to 0 at the root 
$z = i|\xi|$ gives
$$-|\xi|\sum h_{kk} - 2i\xi_{k}h_{0k} = 0.$$
Using again $\sum h_{kk} = tr h - h_{00}$ on the first term and \eqref{e3.10} 
on the second term gives $-|\xi|tr h + |\xi|h_{00} - 2|\xi|h_{00} + |\xi|tr h = 0$, 
which implies that $h_{00} = 0$, and again \eqref{e3.9} then gives $h = 0$. 

  A similar calculation shows that the boundary data \eqref{e3.6}-\eqref{e3.7} 
may be continuously deformed to full Dirichlet boundary data $g_{\alpha\beta} = 
\gamma_{\alpha\beta}$ maintaining ellipticity, cf.~also ~\cite{S}. The latter boundary 
value problem clearly has index 0, and hence, by the homotopy invariance of the 
index, so does the boundary system \eqref{e3.6} or \eqref{e3.7}. 
{\endproof}

  Next we consider some applications of Proposition 3.1. Probably the most 
natural choice for the form $B$ is just $B = g^{T}$, so that for $\gamma = g^{T}$, 
$[\gamma]_{B} = [\gamma]$ is the conformal class of $\gamma$. This leads to 
Theorem 1.2. 

\medskip

{\bf Proof of Theorem 1.2.}

  Let ${\mathcal C}^{m,\alpha}(\partial M)$ be the space of pointwise 
conformally equivalent $C^{m,\alpha}$ metrics on $\partial M$. Proposition 3.1 
and elliptic boundary regularity, cf.~\cite{ADN, Mo}, implies that the map
\begin{equation}\label{e3.13}
\Psi: Met_{C}^{m,\alpha}(M) \rightarrow S_{2}^{m-2,\alpha}(M)\times 
{\mathcal C}^{m,\alpha}(\partial M)\times C^{m-1,\alpha}(\partial M),
\end{equation}
$$\Psi(g) = (\Phi_{\widetilde g}(g), \ [g^{T}], \ H),$$
is a smooth Fredholm map of index 0 for $g$ near $\widetilde g$. Hence, 
the associated boundary map
\begin{equation}\label{e3.14}
\widetilde \Pi_{D}: {\mathbb E}^{m,\alpha}_{C}(M) \rightarrow 
{\mathcal C}^{m,\alpha}(\partial M)\times C^{m-1,\alpha}(\partial M),
\end{equation}
$$\widetilde \Pi_{D}(g) = ([g^{T}], \ H),$$
is also smooth and Fredholm, of Fredholm index 0 for $g$ near $\widetilde g$. 
The proof of Theorem 1.2 then follows from Lemma 2.9, just as in the proof 
of Theorem 1.1. 
{\endproof}

  Note that for $g \in {\mathcal E}^{m,\alpha}$, the scalar constraint \eqref{e3.4} 
implies that the scalar curvature of the boundary metric $\gamma = g^{T} \in 
Met^{m,\alpha}(\partial M)$ is in $C^{m-1,\alpha}$. This is consistent with the 
fact that only the conformal class of the boundary metric $g^{T}$ is prescribed 
in \eqref{e3.14}. 

\medskip

  Next, consider the example where $B$ equals the $2^{\rm nd}$ fundamental form 
$A$ of the metric $g \in {\mathbb E}_{C}^{m,\alpha}$ and assume $\partial M$ 
is strictly convex for $(M, g)$. One has $A \in S_{2}^{m-1,\alpha}(\partial M)$, 
so that the quotient $Met^{m,\alpha}(\partial M)/A$ is not well-defined. To 
remedy this, let $\widetilde A = \widetilde A(g, \varepsilon) \in S_{2}^{\infty}
(\partial M)$ be a $C^{\infty}$ smooth approximation to $A = A(g)$, 
$\varepsilon$-close to $A$ in $C^{m-1,\alpha}$. As above, the 
boundary map
\begin{equation}\label{e3.15}
\widetilde \Pi_{\widetilde A}: {\mathbb E}^{m,\alpha}_{C}(M) \rightarrow 
Met^{m,\alpha}(\partial M)/ \widetilde A \times C^{m-1,\alpha}(\partial M),
\end{equation}
$$\widetilde \Pi_{D}(g) = ([g^{T}]_{\widetilde A}, \ H),$$
is $C^{\infty}$ smooth and Fredholm, of Fredholm index 0 for $g$ near 
$\widetilde g$. In particular, the linearized map has finite dimensional 
kernel and cokernel. This leads to the following result, closely related 
to a result of Schlenker ~\cite{S}. 

\begin{proposition}\label{p3.2}
Suppose $\partial M$ is strictly convex in $(M, g)$. Then near $g$, the space 
of boundary values ${\mathcal B} = \Pi_{D}({\mathcal E}^{\infty})$ of $C^{\infty}$ 
Einstein metrics on $M$, if non-empty, is a variety of finite codimension in 
$Met^{\infty}(\partial M)$. 
\end{proposition}

{\bf Proof:} It suffices to prove the result at the linearized level. 
First, note that the Fredholm property of the boundary map \eqref{e3.15} also 
holds when $m = \infty$. Observe also that the full diffeomorphism group 
${\mathcal D}^{\infty}$ acts on ${\mathbb E}^{\infty}$, (but not on the slice 
${\mathbb E}_{C}^{\infty}$). 

   Suppose then $h \in Im D\Pi \subset S_{2}^{\infty}(\partial M)$. The projection 
$S_{2}^{\infty}(\partial M) \rightarrow S_{2}^{\infty}(\partial M) / \widetilde A$ 
sends $Im D\Pi$ onto a subspace of finite codimension. On the other hand, 
regarding the fiber of this projection, for $f \in C^{\infty}$, one has 
$fA = \delta^{*}(fN) \in S_{2}^{\infty}(\partial M)$, so that 
$h + fA \in Im D\Pi$, for any such $f$ and any $h \in Im D\Pi$. It follows that 
$Im D\Pi$ is $\varepsilon$-dense in a subspace of finite codimension in 
$S_{2}^{\infty}(\partial M)$. One may then let $\widetilde A = 
\widetilde A(\varepsilon) \rightarrow A$ in $C^{\infty}$, and the 
result follows. 
{\endproof}

  We point out that natural analogs of Propositions 3.1, 3.2 and the discussion 
above also for Neumann boundary value problems, (replacing $\gamma$ by $A$). 
The details of this are left to the interested reader.

\section{Extension to complete, noncompact metrics.}
\setcounter{equation}{0}

 In this section, we consider extensions of Theorems 1.1 and 1.2 to 
complete open manifolds with compact boundary. Of course this is only 
relevant in the case $\lambda \leq 0$, since Einstein metrics of positive 
Ricci curvature have a bound on their diameter.  

 Let $M$ be an open manifold with compact boundary, in the sense that 
$M$ has a compact (interior) boundary $\partial M$, together with a 
collection of non-compact ends. Apriori, at this stage $M$ could have 
an infinite number of ends, and/or ends of infinite topological type. 
As in \S 2, we assume $\pi_{1}(M, \partial M) = 0$, so that in particular 
$\partial M$ is connected. 

 Let $g_{0}$ be an Einstein metric on $M$ which is $C^{m,\alpha}$ up to 
$\partial M$, $m \geq 2$, and which is complete away from $\partial M$. 
Choose also a fixed, locally finite atlas in which the metric $g_{0}$ is 
locally in $C^{m,\alpha}$ up to $\partial M$. 

 The metric $g_{0}$ determines the asymptotic behavior of the space of 
metrics to be considered. To describe this, on $(M, g_{0})$, let
\begin{equation} \label{e4.1}
v(r) = vol S(r), 
\end{equation}
where $S(r)$ is the geodesic $r$-sphere about $\partial M$, i.e.~$S(r) = 
\{x \in (M, g_{0}): dist(x, \partial M) = r\}$. Choose positive constants, 
$a, b > 0$ and let $Met_{0}(M) = Met_{g_{0},a,b}^{m,\alpha}(M)$ be the space 
of $C^{m,\alpha}$ metrics on $M$, (in the given atlas), such that, for $r$ large, 
\begin{equation} \label{e4.2}
|g - g_{0}|(r) = \sup_{x\in S(r)}|g - g_{0}|(x) \leq  r^{-a}, 
\end{equation}
\begin{equation} \label{e4.3}
|\nabla^{k}g|(r) = \sup_{x\in S(r)}|\nabla^{k}g|(x) \leq  r^{-(a+b)}, 
\end{equation}
for $k = 1,2$ and for any $g_{1}, g_{2} \in  Met_{0}(M)$,
\begin{equation} \label{e4.4}
|g_{1} - g_{0}|\cdot |\nabla g_{2}|(r) + |\nabla g_{1}|\cdot 
|\nabla^{2}g_{2}|(r)  \leq  \varepsilon (r)v(r)^{-1}, 
\end{equation}
where $\varepsilon (r) \rightarrow  0$ as $r \rightarrow \infty$; the norms 
and covariant derivatives are taken with respect to $g_{0}$. 

 These decay conditions at infinity are quite weak. Consider for 
example the situation where $g_{0}$ is Euclidean, or more generally 
flat in the sense that $g_{0}$ is the flat metric on ${\mathbb R}^{m}
\times T^{n-m+1}$, where $T^{n-m+1}$ is a flat $(n-m+1)$ torus. Then 
$v(r) = cr^{m-1}$ and the conditions (4.2)-(4.4) are satisfied if 
\begin{equation} \label{e4.5}
2a+b > m-1. 
\end{equation}
The usual notion of an asymptotically flat metric $g$ requires $g$ to 
decay at the rate of the Green's function, $|g - g_{0}| = 
O(r^{-(m-2)})$ in this case, while $|\nabla^{k}g| = O(r^{-(m-2+k)})$. 
The condition (4.5) is clearly much weaker than this requirement. 

 Now given $g_{0}$, define the spaces $Met_{C}$, $Z_{C}$ and ${\mathbb E}_{C} 
\subset  Z_{C}$ as subspaces of $Met_{0}(M)$ exactly as in \eqref{e2.9}-\eqref{e2.10}. 
Further, let 
\begin{equation} \label{e4.6}
{\mathbb E} \subset  Met_{0}(M), 
\end{equation}
be the space of all Einstein metrics in $Met_{0}(M)$. 

  Next, regarding the gauge groups for these spaces, let ${\mathcal D}$ be 
the group of $C^{m+1,\alpha}$ diffeomorphisms $\phi$ of $M$ which satisfy 
decay conditions analogous to (4.2)-(4.4), i.e.~taking the supremum over 
$x \in S(r)$,  
\begin{equation} \label{e4.7}
|\phi  - Id|(r) \leq  r^{-a}, 
\end{equation}
\begin{equation} \label{e4.8}
|\nabla^{k}(\phi  - Id)|(r) \leq  r^{-(a+b)}, 
\end{equation}
for $k = 1,2$ and for any $\phi_{1}, \phi_{2} \in  {\mathcal D}$,
\begin{equation} \label{e4.9}
|\phi_{1} - Id|\cdot |\nabla (\phi_{2} - Id)|(r) + |\nabla (\phi_{1} - 
Id)|\cdot |\nabla^{2}(\phi_{2} - Id)|(r) \leq  \varepsilon (r)v(r)^{-1}. 
\end{equation}
Then ${\mathcal D}$ acts on $Met_{0}(M)$. Let ${\mathcal D}_{1} \subset 
{\mathcal D}$ be the subgroup of diffeomorphisms equal to the indentity on 
$\partial M$. Let $\chi_{1}$ denote the corresponding space of vector 
fields on $M$. 

   The proofs of Theorems 1.1-1.2 in this context are identical to the proofs 
when $M$ is compact, provided two issues are addressed. First, in the 
integration by parts arguments used in several places in the proof of 
Proposition 2.5 and the lemmas preceding it, one needs all boundary terms 
taken over $S(r)$ to decay to $0$ as $r \rightarrow \infty$. Second, one 
needs to choose function spaces and boundary conditions at infinity for which 
the operators $L$ and $\beta \delta^{*}$ are Fredholm. 

 Thus, consider closed subspaces $Met_{F}(M) \subset Met_{C}(M)$, 
${\mathcal D}_{F} \subset {\mathcal D}_{1}$, and the associated $S_{F}(M) \subset 
S_{2}(M)$ and $\chi_{F}(M) \subset \chi_{1}$, which are compatible in the 
sense that ${\mathcal D}_{F}$ acts on $Met_{F}(M)$. One may then consider the 
quotient spaces $Met_{F}(M)/{\mathcal D}_{F}$ and in particular
\begin{equation}\label{e4.10}
{\mathcal E}_{F} = {\mathbb E}_{F}/{\mathcal D}_{F},
\end{equation}
where ${\mathbb E}_{F} \subset {\mathbb E}$ is the subspace of Einstein metrics 
in $Met_{F}(M)$. 

\begin{proposition}\label{p4.1} 
Let $Met_{F}(M) \subset  Met_{C}(M)$ and ${\mathcal D}_{F} \subset {\mathcal D}_{1}$ 
be compatible closed subspaces on which the operators $L|_{T_{g}Met_{F}(M)}$ 
and $\beta\delta^{*}|_{\chi_{F}}$ are Fredholm.

  Then Theorems 1.1-1.2 hold on $Met_{F}(M)$, i.e.~the space ${\mathcal E} 
= {\mathcal E}_{F}(M)$ is either empty or an infinite dimensional smooth Banach 
manifold, (or Fr\'echet manifold when $m = \infty$), on which the boundary 
map \eqref{e1.5} satisfies the conclusions of Theorem 1.2. 
\end{proposition}

{\bf Proof:} By straightforward inspection, the decay conditions (4.2)-(4.4) 
and (4.7)-(4.9) insure that the first condition above regarding the decay of 
the boundary terms at $S(r)$ holds. These boundary terms arise in \eqref{e2.11}, 
\eqref{e2.13} and in \eqref{e2.21}, via the divergence theorem. 

  Given that $L$ and $\beta\delta^{*}$ is Fredholm, the proofs of Theorems 
1.1-1.2 then carry over without change to the current situation. 
{\endproof}

 For an arbitrary complete Einstein metric $(M, g_{0})$, there is no 
general theory to determine whether natural elliptic operators are 
Fredholm on suitable function spaces. A detailed analysis in the 
case of ``fibered boundary'' metrics has been carried out by Mazzeo 
and Melrose, cf.~\cite{MM}. For simplicity, we restrict here to the situation 
of asymptotically flat metrics. 

 Thus, let $g_{fl}$ be a complete flat metric on the manifold $N = 
{\mathbb R}^{m}\times T^{n+1-m}/ \Gamma$, where $\Gamma$ is a finite group 
of isometries. Let $(x,y)$ be standard coordinates for ${\mathbb R}^{m}$ and 
$T^{n+1-m}$ and let $r = |x|$. Define
$$C_{\delta}^{m,\alpha}(N) = \{u = r^{-\delta}f: f \in C_{0}^{m,\alpha}(N)\}, $$
where $C_{0}^{m,\alpha}$ is the space of functions $f$ such that 
$(1 + r^{2})^{|\beta|/2}\partial_{x}^{\beta}f \in  C^{0,\alpha}, 
\partial_{y}^{\beta}f \in  C^{0,\alpha}$, where $|\beta| \leq  m$, and 
$C^{0,\alpha}$ is the usual space of $C^{\alpha}$ Holder continuous 
functions on $N$. 

  Let $M$ be a manifold with compact boundary $\partial M$, having a 
finite number of ends, each diffeomorphic to some $N$ above, (not 
necessarily fixed). Given a choice of flat metric $g_{fl}$ on each end, 
let $Met_{\delta}(M)$ be the space of locally $C^{m,\alpha}$ metrics $g$ 
on $M$ such that the components of $(g - g_{fl})$ in the $(x,y)$ coordinates 
are in $C_{\delta}^{m,\alpha}(N)$. One defines the group of $C^{m+1,\alpha}$ 
diffeomorphisms ${\mathcal D}_{1,\delta}$ and associated vector fields 
$\chi_{1,\delta}$ in the same way. 

 By ~\cite{MM}, the Laplace-type operators $L$ and $\delta\delta^{*}$ are 
Fredholm as maps $Met_{\delta}^{m,\alpha}(M)\rightarrow 
S_{2,\delta}^{m-2,\alpha}(M)$ and $\chi_{1,\delta}^{m+1,\alpha} 
\rightarrow \chi_{1,\delta}^{m-1,\alpha}$ provided 
\begin{equation} \label{e4.11}
0 <  \delta  <  m-2.
\end{equation}
Choosing then $a = \delta$ and $b = \delta +1$ in (4.2)-(4.3) and 
(4.7)-(4.8) shows that (4.4) and (4.9) hold provided 
\begin{equation}\label{e4.12}
\frac{m-2}{2} <  \delta  <  m-2.
\end{equation}

  We now set $Met_{F}(M) = Met_{\delta}^{m,\alpha}(M)$, for $\delta$ satisfying 
(4.12) and let ${\mathcal D}_{F}$ be the corresponding space of $C^{m+1,\alpha}$ 
diffeomorphisms. Let ${\mathcal E} = {\mathcal E}_{F}$. Then combining the 
results above with the rest of the proof of Theorem 1.1 proves the following 
more precise version of Theorem 1.3. 
\begin{theorem}\label{t4.2}
For $\pi_{1}(M, \partial M) = 0$, the space ${\mathcal E}$ of Ricci-flat, locally 
asymptotically flat metrics on $M$, satisfying the decay conditions {\rm (4.12)}, 
if non-empty, is an infinite dimensional smooth Banach manifold, (Fr\'echet if 
$m = \infty$). Further, the boundary map \eqref{e1.5} satisfies the conclusions 
of Theorem 1.2. 
\end{theorem}
{\endproof}

  Note that Einstein metrics $g \in {\mathcal E}$ will often satisfy stronger decay 
conditions than (4.12). The Einstein equations imply that the metrics decay to the 
flat metric on the order of $O(r^{-(m-2)})$; this will not be discussed further 
here however.

\section{Matter fields.}
\setcounter{equation}{0}

 In this section, we consider Theorems 1.1 - 1.2 for the Einstein equations 
coupled to other (matter) fields $\phi$. Typical examples of such fields, which 
arise naturally in physics are:

$\bullet$ Scalar fields, $u: M \rightarrow {\mathbb R}$.

$\bullet$ $\sigma$-models, $\varphi: (M, g) \rightarrow (X, \sigma)$, where $(X, \sigma)$ 
is a Riemannian manifold. 

$\bullet$ Gauge fields $A$, i.e. connection 1-forms, with values in a Lie algebra, 
on principal bundles over $M$. 

$\bullet$ $p$-form fields $\omega$. 

  We assume that there is an action or Lagrangian ${\mathcal L} = {\mathcal L}(g, \phi)$, of 
the form 
\begin{equation}\label{e5.1}
{\mathcal L} = {\mathcal L}_{EH} + {\mathcal L}_{m},
\end{equation}
where ${\mathcal L}_{EH}$ is the 
Einstein-Hilbert Lagrangian with integrand $(s - 2\Lambda)dV$ and where the 
matter Lagrangian ${\mathcal L}_{m}$ involves the fields $\phi$ up to $1^{\rm st}$ order, 
with coupling to the metric $g$ also involving at most the $1^{\rm st}$ derivatives 
of $g$. We also assume that ${\mathcal L}$ is analytic in $(g, \phi)$ and is diffeomorphism 
invariant, in that for any $f \in {\mathcal D}_{1}$, 
\begin{equation}\label{e5.2}
{\mathcal L}(f^{*}g, f^{*}\phi) = {\mathcal L}(g, \phi).
\end{equation}

  The variation of ${\mathcal L}$ with respect to $g$, $\frac{\partial {\mathcal L}}
{\partial g}$ gives the Euler-Lagrange equations for $g$:
\begin{equation} \label{e5.3}
-E^{1}(g, \phi) = Ric_{g} - \frac{s}{2}g + \Lambda g - T = 0,
\end{equation}
where $T$ is the stress-energy tensor of the fields $\phi$, i.e.~the variation of 
${\mathcal L}_{m}$ with respect to $g$, cf.~\cite{HE} for instance. The stress-energy 
$T$ is $1^{\rm st}$ order in $g$ and $\phi$ and the Bianchi identity implies the 
conservation property
\begin{equation}\label{e5.4}
\delta T = 0.
\end{equation}
   Similarly, the variation of ${\mathcal L}_{m}$ with respect to the fields $\phi$ gives 
the Euler-Lagrange equations for $\phi$, written schematically as 
\begin{equation} \label{e5.5}
E^{2}(g, \phi) = E_{g}^{2}(\phi) = 0.
\end{equation}
We assume $E_{g}^{2}(\phi)$ can be written in the form of a $2^{\rm nd}$ order elliptic 
system for $\phi$, with coefficients depending on $g$ up to $1^{\rm st}$ order. 
Typically, the operator $E_{g}^{2}$ will be a diagonal or uncoupled system of Laplace-type 
operators at leading order. For simplicity, we do not discuss Dirac-type operators, 
although it can be expected that similar results hold in this case. Note that 
by (5.2), the coupled field equations (5.3) and (5.5) are invariant under the action 
of ${\mathcal D}_{1}$. 

   For example, the Lagrangian for a scalar field with potential $V$ is given by 
\begin{equation} \label{e5.6}
{\mathcal L}_{m} = -\int_{M}[{\tfrac{1}{2}}|du|^{2} + V(u)]dV_{g},
\end{equation}
where $V: {\mathbb R} \rightarrow {\mathbb R}$. An important special case is the 
free massive scalar field, where $V(u) = m^{2}u^{2}$. The field equation (5.5) 
for $u$ is then
\begin{equation} \label{e5.7}
\Delta_{g}u  = V'(u),
\end{equation}
with stress-energy tensor given by
\begin{equation} \label{e5.8}
T = {\tfrac{1}{2}}[du \cdot du - ({\tfrac{1}{2}}|du|^{2} + V(u))g] .
\end{equation}

   For a gauge field or connection 1-form $\phi = d + A$, the usual Lagrangian is 
the Yang-Mills action 
\begin{equation} \label{e5.9}
{\mathcal L}_{m} = -{\tfrac{1}{2}}\int_{M}|F|^{2}dV_{g},
\end{equation}
where $F = d_{A}A \equiv dA + \frac{1}{2}[A, A]$ is the curvature of $A$. The field equations 
are the Yang-Mills equations, (or Maxwell equations in the case of a $U(1)$ bundle):
\begin{equation} \label{e5.10}
d_{A}F = \delta_{A}F = 0,
\end{equation}
with stress-energy tensor
\begin{equation}\label{e5.11}
T = F \cdot F - {\tfrac{1}{2}}|F|^{2}g,
\end{equation}
where $(F \cdot F)_{\mu\nu} = \langle F_{\mu\alpha}, F_{\nu\beta} \rangle g^{\alpha\beta}$. 

  To match with the work in \S 2, we pass from (5.3) to the equivalent equations
\begin{equation}\label{e5.12}
Ric_{g} - \lambda g - T_{0} = 0,
\end{equation}
where $T_{0} = T - \frac{tr T}{n+1}g$ is the trace-free part of $T$. The conservation 
law (5.4) then translates to 
\begin{equation}\label{e5.13}
\beta(\widetilde T) = 0.
\end{equation}

   We begin with a detailed discussion of the case of the Einstein equations 
coupled to a scalar field $u: M \rightarrow {\mathbb R}$ with potential $V(u)$, 
where $V: {\mathbb R} \rightarrow {\mathbb R}$ is an arbitrary smooth function; 
as will be seen below, the treatment of other fields is very similar. 
 
  The full Lagrangian is given by 
\begin{equation}\label{e5.14}
{\mathcal L} (g, u) = \int_{M}[(s - 2\Lambda) - {\tfrac{1}{2}}|du|^{2} - V(u)]dV_{g},
\end{equation}
which gives the field equations 
\begin{equation}\label{e5.15}
Ric_{g} - \lambda g = T_{0} = {\tfrac{1}{2}}(du \cdot du - {\tfrac{1}{n-1}}Vg), 
\ \ \Delta u = V'(u),
\end{equation}
when the variations of $(g, u)$ are of compact support in $M$. As in the proof of 
Theorem 1.1, where the boundary data for the metric $g$ were not fixed in advance, 
it is useful here not to fix boundary values for the scalar field $u$. Thus, instead 
of (5.14), we consider the Lagrangian
\begin{equation}\label{e5.16}
{\mathcal L} (g, u) = \int_{M}[(s - 2\Lambda) + {\tfrac{1}{2}}u\Delta u - V(u)]dV.
\end{equation}
Of course, the Lagrangians (5.14) and (5.16) differ just by boundary terms. 

  The Lagrangian is a map $Met^{m,\alpha}(M)\times C^{k,\beta}(M) \rightarrow {\mathbb R}$, 
and we assume $k \geq 2$, $\beta \in (0,1)$. The differential (or variation) $d{\mathcal L}$ 
is then a map 
\begin{equation}\label{e5.17}
d{\mathcal L}  = ({\mathcal L}^{1}, {\mathcal L}^{2}): Met^{m,\alpha}(M)\times C^{k,\beta}(M) 
\rightarrow  T^{*}(Met^{m-2,\alpha}(M)\times C^{k-2,\beta}(M)),
\end{equation}
where
\begin{equation}\label{e5.18}
d{\mathcal L}_{(g,u)}^{1}(h,v) = -\langle Ric_{g} - \lambda g - \widetilde T(g,u), h\rangle dV,
\end{equation}
represents the variation with respect to $g$ and 
\begin{equation}\label{e5.19}
d{\mathcal L}_{(g,u)}^{2}(h,v) =  ({\tfrac{1}{2}}(v\Delta u + u\Delta v) - 
V'(u)\cdot v)dV,
\end{equation}
represents the variation with respect to $u$. The ``new'' stress-energy tensor 
$\hat T$ for (5.16) is given by $\hat T = \frac{1}{2}[u\Delta' u + (\frac{1}{2}u\Delta u - 
V(u))g]$, where $\Delta' u$ is the metric variation of the Laplacian, 
given by 
\begin{equation}\label{e5.20}
\Delta' u(h) = -\langle D^{2}u, h\rangle  + \langle du, \beta (h)\rangle ,
\end{equation}
where $\beta$ is the Bianchi operator $\beta (h) = \delta h + \frac{1}{2}d(tr h)$. 
Using (5.20), it is easily seen that $\hat T = T$, for $T$ as in (5.8), modulo 
boundary terms, and so we continue to use (5.8). In particular, for variations of 
compact support, one obtains the Euler-Lagrange equations (5.15). 

\medskip

 Let ${\mathbb E} = {\mathbb E}(M, g, u)$ denote the space of all solutions of 
the equation $d{\mathcal L} = 0$, for ${\mathcal L}$ as in (5.16), i.e.~the 
space of solutions to the Einstein equations coupled to the scalar field $u$. 
This space is invariant under the diffeomorphism group ${\mathcal D}_{1}$, 
acting on both $(g, u)$ by pullback. 

 As in \S 2, one needs to choose a gauge to break the diffeomorphism 
invariance; (the scalar field has no internal symmetry group, so there 
is no need of an extra gauge for $u$). Thus, analogous to the 
discussion in \S 2, given a background metric $\widetilde g \in {\mathbb E}$, 
define
\begin{equation}\label{e5.21}
\Phi  = \Phi_{\widetilde g}: Met (M)\times C^{k,\beta}(M) \rightarrow  
T^{*}(Met (M)\times C^{k-2,\beta}(M));
\end{equation}
$$\Phi (g, u) = [(Ric_{g} - \lambda g - \widetilde T(g,u) + 
\delta^{*}\beta_{\widetilde g}(g))dV, -({\tfrac{1}{2}}(\cdot \Delta u + 
u\Delta \cdot ) - V'(u)\cdot )dV].$$
(For convenience, we have switched the signs in comparison with (5.18)-(5.19)). 
As before, $Met_{C}(M) = Met_{C}^{m,\alpha}(M)$ is defined to be the space of 
$C^{m,\alpha}$ metrics satisfying the Bianchi constraint \eqref{e2.9}, and we set
$$Z_{C} = \Phi^{-1}(0) \subset  Met_{C}(M)\times C^{m,\alpha}(M).$$
Corollary 2.3 also holds as before, so that 
\begin{equation}\label{e5.22}
Z_{C} = {\mathbb E}_{C},
\end{equation}
(for any boundary conditions on $u$). Given this, the main task is to verify 
that the analog of Proposition 2.5 holds. 

\begin{proposition}\label{p5.1}
Proposition {\rm 2.5} holds for the map $\Phi$ in {\rm (5.21)}, i.e.~$D\Phi$ is 
surjective. 
\end{proposition}

{\bf Proof:} Consider the derivative of $\Phi$ at $g = \widetilde g$:
$$D\Phi  = (D\Phi^{1}, D\Phi^{2}): T(Met(M)\times C^{m,\alpha}(M)) 
\rightarrow  T(T^{*}(Met(M)\times C^{m-2,\alpha}(M))). $$
This is a block matrix of the form 
\begin{equation}\label{e5.23}
{\mathcal H} = 
\left(
\begin{array}{ccc}
\frac{\partial\Phi^{1}}{\partial g} &  & \frac{\partial\Phi^{1}}{\partial u} \\
 &  & \\
\frac{\partial\Phi^{2}}{\partial g} &  & \frac{\partial\Phi^{2}}{\partial u} 
\end{array}
\right)
\end{equation}
The matrix ${\mathcal H}$ is essentially the same as the $2^{\rm nd}$ variation 
of the Lagrangian (5.16); they agree modulo the gauge term $\delta^{*}
\beta_{\widetilde g}(g)$. A straightforward computation, using the fact 
that $(g, u) \in {\mathbb E}$, gives: 
\begin{equation}\label{e5.24}
\frac{\partial\Phi^{1}}{\partial g}(h) = \widetilde L(h) = L(h) + S(h) ,
\end{equation}
\begin{equation}\label{e5.25}
\frac{\partial\Phi^{1}}{\partial u}(v) = -du\cdot dv + 
{\tfrac{1}{2}}V'(u)v g,
\end{equation}
\begin{equation}\label{e5.26}
\frac{\partial\Phi^{2}}{\partial g}(h) = - {\tfrac{1}{2}}(\cdot \Delta' u + 
u\Delta'\cdot  ) - [{\tfrac{1}{4}}(\cdot \Delta u + u\Delta\cdot ) - 
{\tfrac{1}{2}}V'(u) \cdot ]trh
\end{equation}
\begin{equation}\label{e5.27}
\frac{\partial\Phi^{2}}{\partial u}(v) = -{\tfrac{1}{2}}(\cdot \Delta v + 
v\Delta\cdot ) + V''(u)v\cdot  ,
\end{equation}
where $L$ is the Bianchi-gauged linearized Einstein operator \eqref{e2.6} and $S(h)$ is 
an algebraic operator of the form
$$S(h) = {\tfrac{1}{2}}tr hdu\cdot du - {\tfrac{1}{2}}V(u)h - {\tfrac{1}{2}}tr hV(u)g.$$

 Now if there exists $(k, w) \perp Im(D\Phi)$, then 
\begin{equation}\label{e5.28}
\int_{M}\langle D\Phi (h, v), (k,w)\rangle dV = 0,
\end{equation}
for all $(h,v)$ with $h\in T(Met(M)\times C^{k,\beta}(M))$. As in the proof 
of Proposition 2.5, one integrates the expressions (5.24)-(5.27) by parts. For 
(5.24), one obtains, as in \eqref{e2.21}, 
\begin{equation}\label{e5.29}
\int_{M}\langle L(h), k\rangle  + \langle S(h), k \rangle = \int_{M}\langle L(k), h\rangle + 
\langle S(k), h \rangle + \int_{\partial M}D(h, k),
\end{equation}
where $D$ is given by \eqref{e2.22}. For (5.25):
$$[-du\cdot dv + V'(u) \cdot v)g](k) = 
\int_{M}-\langle du \cdot dv, k\rangle + {\tfrac{1}{2}}V'(u)\cdot v trk =  $$
\begin{equation}\label{e5.30}
\int_{M}v[-\delta (k(du)) + {\tfrac{1}{2}}V'(u) trk] -\int_{\partial M}vk(du,N).
\end{equation}
Next for (5.26):
$$-[{\tfrac{1}{2}}(\cdot \Delta' u + u\Delta'\cdot) + 
{\tfrac{1}{4}}(\cdot \Delta u + u\Delta\cdot ) - {\tfrac{1}{2}}V'(u)\cdot )trh](w) = $$
$$-\int_{M}  {\tfrac{1}{2}}(w\Delta' u + u\Delta' w) + {\tfrac{1}{4}}
(w\Delta u + u\Delta w) - {\tfrac{1}{2}}V'(u)w)trh $$
and
$$\int_{M}w\Delta' u = \int_{M}-w\langle D^{2}u, h\rangle + w\langle du, \beta (h)\rangle = $$
$$\int_{M}\langle du \cdot dw + {\tfrac{1}{2}}\delta (wdu)g, h\rangle - 
\int_{\partial M}wh(du,N) - {\tfrac{1}{2}}trh wN(u).$$
Interchanging $u$ and $w$ then gives
\begin{equation}\label{e5.31}
w\cdot \frac{\partial\Phi^{2}}{\partial g}(h) = -\int_{M}\langle du \cdot dw + 
{\tfrac{1}{4}}\delta (duw)g + {\tfrac{1}{4}}(w\Delta u + u\Delta w)g - 
{\tfrac{1}{2}}V'(u)w)g, h\rangle
\end{equation}
$$+ {\tfrac{1}{2}}\int_{\partial M}h(duw, N) - {\tfrac{1}{2}}trh N(uw).$$
Finally, for (5.27), 
$$-[{\tfrac{1}{2}}(\cdot \Delta v + v\Delta\cdot ) - V''(u)v\cdot ](w) = 
-{\tfrac{1}{2}}\int_{M}(w\Delta v + v\Delta w - 2V''(u)vw)$$
\begin{equation}\label{e5.32}
= -\int_{M}v(\Delta w - V''(u)w) - {\tfrac{1}{2}}\int_{\partial M}wN(v) - vN(w).
\end{equation}

 Now, supposing (5.28) holds, since $v$ is arbitrary, by adding the bulk terms in 
(5.30) and (5.32) one obtains
\begin{equation}\label{e5.33}
\Delta w - V''(u)\cdot w + \delta (k(du)) - {\tfrac{1}{2}}V'(u)trk = 0,
\end{equation}
on $(M, g)$; this is the equation for the variation $w$ of the scalar 
field $u$. Adding the boundary terms in (5.30) and (5.32) gives
\begin{equation}\label{e5.34}
\int_{\partial M}wN(v) - v[N(w) - 2k(du,N)] = 0.
\end{equation}
The boundary values of $v$ are arbitrary, so that both $v$ and 
$N(v)$ can be prescribed arbitrarily at $\partial M$. Hence (5.34) implies that
\begin{equation}\label{e5.35}
w = 0 \ \ {\rm and} \ \  N(w) - 2k(du,N) = 0,
\end{equation}
at $\partial M$. 

  Next, since $h$ is arbitrary in the interior, adding the bulk terms in (5.29) 
and (5.31) gives
\begin{equation}\label{e5.36}
L(k) + S(k) - du \cdot dw - {\tfrac{1}{4}}\delta (duw)g - {\tfrac{1}{4}}(w\Delta u - 
u\Delta w)g + {\tfrac{1}{2}}V'(u)wg = 0.
\end{equation}
This is the equation for the variation $k$ of the metric $g$. At $\partial M$, 
adding the boundary terms in (5.29) and (5.31) gives 
$$D(h,k) + {\tfrac{1}{2}}[h(duw,N) - {\tfrac{1}{2}}trh N(uw)] = 0. $$
Since $w = 0$ at $\partial M$, one thus has
\begin{equation}\label{e5.37}
D(h,k) + {\tfrac{1}{2}}uN(w)[h(N,N) - {\tfrac{1}{2}}trh] = 0.
\end{equation}

   Now the same arguments as in \eqref{e2.21}-\eqref{e2.37} carry over to this 
situation essentially unchanged. The proof of \eqref{e2.24} follows in the same 
way as before, via the diffeomorphism invariance of ${\mathbb E}$. It follows 
then from (5.35) and (5.37) that the geometric Cauchy data vanish at $\partial M$, 
i.e.
\begin{equation}\label{e5.38}
k^{T} = (A_{k}')^{T} = 0, \ {\rm and} \ \ w = N(w) = 0, \ \ {\rm at} \ \ \partial M.
\end{equation}

 By (5.33) and (5.36), the pair $(k, w)$ satisfy the coupled system of 
equations:
\begin{equation}\label{e5.39}
L(k) + S(k) - du \cdot dw - {\tfrac{1}{4}}\delta (duw)g - {\tfrac{1}{4}}(w\Delta u - 
u\Delta w)g + {\tfrac{1}{2}}V'(u)wg = 0,
\end{equation}
\begin{equation}\label{e5.40}
\Delta w - V''(u)\cdot w + \delta (k(du)) - {\tfrac{1}{2}}V'(u)trk = 0.
\end{equation}
Here, $(g,u)$ are fixed, and viewed as (smooth) coefficients, while $(k,w)$ 
are the unknowns. The equations (5.39)-(5.40) express the fact that $(k,w) 
\in  TZ$. Since $k$ is transverse-traceless so that $\beta (k) = 0$, the pair 
$(k, w)$ satisfy the linearized Einstein equations coupled to a scalar field. 

  The unique continuation property, Proposition 2.4, also holds for these linearized 
Einstein equations, since the scalar field $u$ modifies the Einstein equations only 
at first order. Given the vanishing of the geometric Cauchy data in (5.38), the 
proof that 
\begin{equation}\label{5.41}
k = w = 0 \ \ {\rm on} \ \ M,
\end{equation}
proceeds just as before. This proves the surjectivity of $D\Phi$, and the 
proof that the kernel splits is again the same. 
{\endproof}

  Let ${\mathcal E} = {\mathcal E}^{m,\alpha}_{\lambda, V}(g, u)$ be the moduli 
space of Einstein metrics $g$ coupled to a scalar field $u$ with potential $V$ on 
$(M, \partial M)$. As before, one has a natural Dirichlet boundary map $\Pi_{D}$, 
giving Dirichlet boundary values to $u$, or its mixed version $\widetilde \Pi_{D}$ 
as in \eqref{e1.5}. Given Proposition 5.1 and the remarks above, the rest of the work 
in \S 2 and \S 3 carries over unchanged, and proves:
\begin{corollary}\label{c5.2}
Suppose $\pi_{1}(M, \partial M) = 0$. Then the space ${\mathcal E}$ of solutions 
to the Einstein equations coupled to a scalar field with potential $V$, if non-empty, 
is an infinite dimensional smooth Banach manifold, (Fr\'echet when $m = \infty$), for 
which the boundary map $\widetilde \Pi_{D}$ is smooth and Fredholm of index 0, 
i.e.~Theorems 1.1-1.2 hold. 
\end{corollary}
{\endproof}

   Similarly, Corollary 5.2 holds in the same way, for scalar fields in the space 
$C_{\delta}^{k,\beta}$ with $\delta$ satisfying (4.12). 

  Next consider the situation of the Einstein equations coupled to a nonlinear 
$\sigma$-model. In this case the field $\phi$ is a smooth function $u: (M, g) \rightarrow 
(X, \sigma)$, with matter Lagrangian
\begin{equation}\label{e5.42}
{\mathcal L}_{m} = -\int_{M}[{\tfrac{1}{2}}|du|^{2} + V(u)]dV_{g},
\end{equation}
where $|du|^{2} = \sigma(du(e_{i}), du(e_{i}))$, for a local orthonormal basis 
$e_{i}$ of $(M, g)$; $du$ is the derivative map of $u$, and $V: X \rightarrow {\mathbb R}$ 
the potential function. 

  The analysis in this case is essentially the same as that of a single scalar field 
discussed above. Probably the simplest way to see this is to isometrically embed 
$(X, \sigma)$, via the Nash embedding theorem, into a large Euclidean space 
${\mathbb R}^{N}$. Then $u: (M, g) \rightarrow (X, \sigma)$ is a vector-valued function 
$u = \{u^{i}\}:(M, g) \rightarrow {\mathbb R}^{N}$, $1 \leq i \leq N$, with the constraint 
that $Im \, u \subset X \subset {\mathbb R}^{N}$. The metric $\sigma$ on $X$ is then just 
the restriction of the Euclidean dot-product metric to $TX$. 

   The $\sigma$-model field equations for the Lagrangian (5.42) are
\begin{equation}\label{e5.43}
\Delta^{T}u = V'(u) = u^{*}(\nabla V),
\end{equation}
where $\Delta^{T}$ is the projection of the Laplacian $\Delta = \Delta_{(M, g)}$, 
acting on the components $u^{i}$ of $u$, onto $TX$. If $S$ denotes the $2^{\rm nd}$ 
fundamental form of $X$ in ${\mathbb R}^{N}$, then (5.43) is equivalent to the system
\begin{equation}\label{e5.44}
\Delta u = S(du,du) + V'(u).
\end{equation}
The stress-energy tensor $T$ has exactly the same form as in (5.8), where 
$du \cdot du$ is the symmetric bilinear form on $M$ given by taking the Euclidean 
dot product of the vector $u = \{u^{i}\}$.

   Given this, it is now straightforward to see that all the computations carried out 
in the case of a single scalar field $u$ carry over without significant change to the 
present constrained, vector-valued field $u$ to give:

\begin{corollary}\label{c5.3}
If $\pi_{1}(M, \partial M) = 0$, then the space ${\mathcal E}$ of solutions to the 
Einstein equations coupled to a $\sigma$-model $u:M \rightarrow (X, \sigma)$, if 
non-empty, is an infinite dimensional smooth Banach manifold, for which the boundary 
map $\widetilde \Pi_{D}$ is smooth and Fredholm, of index 0. 
\end{corollary}
{\endproof}

   Finally consider the Einstein equations coupled to gauge fields, i.e.~connections 
$\omega$ on principal bundles $P$ over $M$ with compact semi-simple structure group $G$ 
with bi-invariant metric. The simplest coupled Lagrangian is
\begin{equation}\label{e5.45}
{\mathcal L} = \int_{M}(s - 2\Lambda)dV_{g} - {\tfrac{1}{2}}\int_{M}|F|^{2}dV_{g},
\end{equation}
with field equations
\begin{equation}\label{e5.46}
Ric - \lambda g - T_{0} = 0, \ \ \delta_{\omega}F = 0,
\end{equation}
where $T$ is given by (5.11) and $T_{0}$ is the trace-free part. 

  Let ${\mathcal A}(P) = {\mathcal A}^{k,\beta}(P)$ denote the space of connections 
on $P$ which are $C^{k,\beta}$ smooth up to $\partial M$, with $k \geq 2$, $\beta 
\in (0,1)$. Given any fixed connection $\omega_{0} \in {\mathcal A}(P)$, any 
$\omega \in {\mathcal A}(P)$ has the form $\omega = \omega_{0} + A$, where $A$ 
is a 1-form on $P$ with values in the Lie algebra ${\mathcal L}(G)$. Let 
${\mathbb E} = {\mathbb E}(g, A)$ be the space of all solutions to the field 
equations (5.46), i.e.~the space of all solutions of the Einstein equations coupled 
to the gauge field $A$. The Lagrangian (5.45) and the field equations (5.46) are 
invariant under the diffeomorphisms ${\mathcal D}_{1}$ of $M$, as well as gauge 
transformations of $P$, again equal to the identity on $\partial M$. We expect 
the natural analogs of Theorems 1.1 - 1.2, (and Theorem 4.2), hold in this 
context as well, by the same methods. However, this will not be discussed 
here in detail, cf.~\cite{Ma} for some discussion along these lines.

\begin{remark}\label{r5.5}
{\rm Although the focus of this work has been on Einstein metrics, the main results 
also apply to other field equations, with the background manifold and metric $(M, g)$ 
arbitrary, (not necessarily Einstein), but fixed. Thus for example, the proof of Corollary 
5.2 shows that the space of solutions to the scalar field equation (5.7) with fixed 
$(M, g)$ is an infinite dimensional smooth Banach manifold, (if non-empty), with 
Dirichlet and Neumann boundary maps Fredholm of index 0. This follows just by 
considering the piece $\partial \Phi^{2} / \partial u$ in ${\mathcal H}$ in (5.23). 
Thus, one may set $k = 0$ following (5.38) and argue as before. Of course in the case 
the potential $V(u)$ is linear, the space of solutions of (5.7) is a linear space. 

   Similarly, the space of harmonic maps $u: (M, g) \rightarrow (X, \sigma)$ with fixed 
data $(M, g)$ and $(X, \sigma)$ also satisfies the conclusions of Theorem 1.1-1.2. This 
has previously been known, cf.~\cite{EL}, only in the case of ``non-degenerate'' harmonic 
maps. }
\end{remark}

\bibliographystyle{plain}

\begin{center}
March, 2008
\end{center}

\smallskip
\noindent
\address{Department of Mathematics\\
S.U.N.Y. at Stony Brook\\
Stony Brook, N.Y. 11794-3651}

\noindent
E-mail: anderson@math.sunysb.edu

\end{document}